\documentclass[11pt]{amsart}
\usepackage{mathrsfs}
\usepackage{amssymb}
\usepackage[all,pdf]{xy}
\usepackage{titletoc}
\usepackage{amscd}
\usepackage{amsmath, amssymb}
\usepackage{amsfonts}
\usepackage[colorlinks,linkcolor=blue,citecolor=blue, pdfstartview=FitH]{hyperref}

\numberwithin{equation}{section}

\theoremstyle{plain}
\newtheorem{thm}{Theorem}[section]
\newtheorem{theorem}[thm]{Theorem}
\newtheorem{lemma}[thm]{Lemma}
\newtheorem{corollary}[thm]{Corollary}
\newtheorem{proposition}[thm]{Proposition}

\theoremstyle{definition}

\newtheorem{remark}[thm]{Remark}

\newtheorem{definition}[thm]{Definition}

\newtheorem{example}[thm]{Example}

\newtheorem{defn-thm}[thm]{Definition-Theorem}

\newcommand{\sA}{{\mathcal A}}
\newcommand{\sB}{{\mathcal B}}
\newcommand{\sC}{{\mathcal C}}
\newcommand{\sD}{{\mathcal D}}

\newcommand{\sH}{{\mathcal H}}
\newcommand{\sI}{{\mathcal I}}

\newcommand{\sO}{{\mathcal O}}

\newcommand{\sR}{{\mathcal R}}

\newcommand{\sT}{{\mathcal T}}

\newcommand{\ssg}{{\mathfrak g}}
\newcommand{\ssn}{{\mathfrak n}}
\newcommand{\ssh}{{\mathfrak h}}

\newcommand{\ssL}{{\mathfrak L}}

\newcommand{\C}{{\mathbb C}}

\newcommand{\F}{{\mathbb F}}

\newcommand{\N}{{\mathbb N}}

\newcommand{\Q}{{\mathbb Q}}

\newcommand{\Z}{{\mathbb Z}}

\newcommand{\Hom}{{\mathrm{Hom}}}
\newcommand{\Ker}{{\mathrm{Ker}}}

\newcommand{\btheorem}{\begin{theorem}}
\newcommand{\etheorem}{\end{theorem}}
\newcommand{\bproposition}{\begin{proposition}}
\newcommand{\eproposition}{\end{proposition}}
\newcommand{\bdefinition}{\begin{definition}}
\newcommand{\edefinition}{\end{definition}}
\newcommand{\bcorollary}{\begin{corollary}}
\newcommand{\ecorollary}{\end{corollary}}
\newcommand{\bproof}{\begin{proof}}
\newcommand{\eproof}{\end{proof}}
\newcommand{\bremark}{\begin{remark}}
\newcommand{\eremark}{\end{remark}}
\newcommand{\eexample}{\end{example}}
\newcommand{\bexample}{\begin{example}}

\newcommand{\elemma}{\end{lemma}}
\newcommand{\blemma}{\begin{lemma}}

\newcommand{\ee}{\end{eqnarray*}}
\newcommand{\be}{\begin{eqnarray*}}

\newcommand{\beq}{\begin{equation}}
\newcommand{\eeq}{\end{equation}}

\newcommand{\lrw}{\longrightarrow}

\usepackage{fancyhdr}
\usepackage{tikz-cd}
\usepackage[all]{xy}
\begin{document}

\title{Applications of spherical twist functors to Lie algebras associated to root categories of preprojective algebras}
\author{Fan Xu, Fang Yang*}
\address{Department of Mathematical Sciences\\
Tsinghua University\\
Beijing 100084, P.~R.~China}
\email{fanxu@mail.tsinghua.edu.cn(F. Xu)}
\address{Department of Mathematical Sciences\\
Tsinghua University\\
Beijing 100084, P.~R.~China}
\email{yangfang19@mails.tsinghua.edu.cn(F. Yang)}

\subjclass[2010]{ 
18E30, 16E35, 16E45, 17B99}
%
\keywords{ 
Preprojective algebra; Ringel-Hall Lie algebra; Weyl group; Spherical
twist functor.
}
\thanks{$*$~Corresponding author.}


\begin{abstract}
   Let $\varLambda_Q$ be the preprojective algebra of a finite acyclic quiver $Q$ of non-Dynkin type and $D^b(\mathrm{rep}^n \varLambda_Q)$ be the bounded derived category of finite dimensional nilpotent $\varLambda_Q$-modules. We define spherical twist functors over  the root category $\sR_{\varLambda_Q}$ of $D^b(\mathrm{rep}^n \varLambda_Q)$ and then realize the Weyl group associated to $Q$ as certain subquotient of the automorphism group of the Ringel-Hall Lie algebra $\ssg(\sR_{\varLambda_Q})$ of $\sR_{\varLambda_Q}$ induced by spherical twist functors. We also present a conjectural relation between certain Lie subalgebras of $\ssg(\sR_{\varLambda_Q})$ and  $\ssg(\sR_Q)$, where $\ssg(\sR_Q)$ is the Ringel–Hall Lie algebra associated to the root category $\sR_Q$ of $Q$.
\end{abstract}

\maketitle

\section{Introduction}
Let $Q$ be a finite (connected) acyclic quiver of non-Dynkin type and $\varLambda_Q$ the preprojective algebra. The preprojective algebra $\varLambda_Q$ plays important roles in studying topics related to the Kac-Moody Lie algebra  $\ssg_Q$ of $Q$ and its enveloping algebra.  Firstly, Lusztig (\cite{Lusztig1991, Lusztig2000}) gave a geometric realization of the positive part $U^+_Q$ of the enveloping algebra of the Kac-Moody Lie algebra $\ssg_Q$ and then constructed a basis of $U^+_Q$ called semicanonical basis indexed by irreducible components of the variety $\varLambda_{\boldsymbol{v}}$ of nilpotent $\varLambda_Q$-modules with dimension vector $\boldsymbol{v}$. Secondly, the preprojective algebra $\varLambda_Q$ is independent of orientations of the quiver $Q$. Thirdly, the preprojective algebra $\varLambda_Q$ is derived 2-Calabi-Yau, i.e. the derived category $D^b(\mathrm{rep}^n\varLambda_Q)$ is a 2-Calabi-Yau category (\cite{Geiss2007},\cite{Keller2008}). Geiss, Leclerc and Schroeer (\cite{Geiss2007}) applied this property to prove the cluster multiplication between evaluation forms. In \cite{Shiraishi2016a}, Shiraishi, Takahashi, and Wada defined the spherical twist functor over $D^b(\mathrm{rep}^n\varLambda_Q)$ and showed  the spherical twist functor can be identified with simple reflections in the  corresponding Weyl group  on the level of the Grothendieck group $K_0(D^b(\varLambda_Q))$ based on the property of derived $2$-Calabi-Yau.

The notion of reflection functors for the categories of representations of quivers was introduced by Bernstein, Gelfand, and Ponomarev \cite{Bernstein1973} whose aim was to obtain a simple and elegant proof of Gabriel's theorem. Applying the BGP-reflection functors in the root categories to the Lie algebra model, Xiao, Zhang and Zhu \cite{Xiao2005} obtained the well-known Weyl group action on the Kac-Moody Lie algebras by explicit formula. However, this construction does not work in general since BGP-reflection functors can only be defined for sinks or sources of quivers. Deng, Ruan and Xiao in \cite{Deng2020} generalized this construction to arbitrary vertices for the star-shaped quivers associated with weighted projective lines. However, they can not give a realization of the composition of two simple reflections, because BGP-reflection functors and mutation functors given in \cite{Deng2020} are functors between different categories.

The aim of this paper is to realize simple reflections in the Weyl group $W_Q$ of the Kac-Moody Lie algebra $\ssg_Q$ as automorphisms of a certain Lie algebra, and then realize the Weyl group $W_Q$ as the subquotient of the automorphism group of the corresponding Lie algebra. The idea is illuminated by \cite{Shiraishi2016a},  where each simple reflection of $W_{Q}$ can be lifted to a spherical twist functor of the bounded derived category $D^b(\mathrm{rep}^n\varLambda_Q)$. 

The paper is organized as follows. In Section 2, we introduce some preliminary results, including the Ringel-Hall Lie algebra, the triangulated hull of an orbit category, etc. As we know, $D^b(\mathrm{rep}^n\varLambda_Q)$ is a 2-Calabi-Yau category, so we can define a spherical twist functor $T_{S_i}$ associated to each simple $\varLambda_Q$-module $S_i$, which induces an autoequivalence on the orbit category $D^b(\mathrm{rep}^n\varLambda_Q)/[2]$. Note that $D^b(\mathrm{rep}^n\varLambda_Q)/[2]$ may not be a triangulated category in general,  using Keller's construction of triangulated hull of orbit categories \cite{Keller2005}, we obtain the triangulated hull $\sR_{\varLambda_Q}$ of the orbit category $D^b(\mathrm{rep}^n\varLambda_Q)/[2]$. In Section 3, we extend the autoequivalence induced by a spherical twist funtor $T_{S_i}$ to a triangle autoequivalence $\bar{T}_{2,S_i}$ of the triangulated hull $\sR_{\varLambda_Q}$, which we still call the spherical twist functor. Similar to \cite{Shiraishi2016a}, the spherical twist functor $\bar{T}_{2,S_i}$ is just the simple reflections associated to the simple root $\alpha_i$ on the Grothendieck group level. On the other hand,  Peng and Xiao (\cite{Peng2000}) constructed the Ringel-Hall Lie algebras associated to $2$-periodic triangulated categories. Applying this construction to the triangulated hull $\sR_{\varLambda_Q}$,  we obtain a Ringel–Hall Lie algebra $\ssg(\sR_{\varLambda_Q})$, even the root category is not proper anymore \cite{Xiao2010}.  We show that the spherical twist functor $\bar{T}_{2,S_i}$ gives rise to an automorphism $\Phi_i$ of $\ssg(\sR_{\varLambda_Q})$ in Section $3$. Thus the Weyl group $W_Q$ is a quotient of the subgroup $\mathrm{Aut}_0(\ssg(\sR_{\varLambda_Q}))$, where $\mathrm{Aut}_0(\ssg(\sR_{\varLambda_Q}))$ is the subgroup of the automorphism group of $\ssg(\sR_{\varLambda_Q})$ generated by $\Phi_i$. In Section 4, we explore connections between the Lie subalgebra $\ssg_0(\sR_Q)$ of $\ssg(\sR_Q)$ and the Lie subalgebra $\ssg_0(\sR_{\varLambda_Q})$ of $\ssg(\sR_{\varLambda_Q})$ generated by simple objects $\hat{u}_{S_i}$, $\hat{u}_{S_i[1]}$ and Cartan elements $\frac{h_{S_i}}{d(S_i)}$ for $i\in Q_0$. We show that the positive (resp. negative) part $\ssn^+_0(\sR_Q)$ (resp. $\ssn^-_0(\sR_Q)$) is isomorphic to the positive (resp. negative) part $\ssn_0^+(\sR_{\varLambda_Q})$ (resp. $\ssn^-_0(\sR_{\varLambda_Q})$). Since the triangulated hull $\sR_{\varLambda_Q}$ is not proper, there are more elements in $\ssg_0(\sR_{\varLambda_Q})$ than in $\ssg_0(\sR_Q)$. Let $\sI$ be the ideal generated by elements $\hat{u}_{E_i}-\hat{u}_{E_i[1]}$, $i\in Q_0$, and $\ssh'\subset \ssh$ be the center of $\ssg_0(\sR_{Q})$. Then $\ssg_0(\sR_{\varLambda_Q})/(\sI+\ssh')$ is either $0$ or  isomorphic to $\ssg_0(\sR_Q)/\ssh'$. 

\section{Preliminaries}
\subsection{Root category}
 Let $\varLambda_Q$ be the preprojective algebra of an acyclic quiver $Q=(Q_0,Q_1)$ over a field $k$,  where $Q$ is not of Dynkin type, then $\varLambda_Q$ is of infinite dimensional with finite global dimension equal to $2$. Let $\mathrm{mod}\varLambda_Q$ be the category of finite dimensional $\varLambda_Q$-modules, and $D^b(\varLambda_Q)$ the bounded derived category of $\mathrm{mod}\varLambda_Q$. Denote by $[1]$ the shift functor. Define the orbit category $D^b(\mathrm{rep}^n\varLambda_Q)/[2]$ of $D^b(\varLambda_Q)$ to be the category having the same objects with $D^b(\mathrm{rep}^n\varLambda_Q)$, and morphism spaces are given by
   $$\Hom_{D^b( \varLambda_Q)/[2]}(X,Y):=\bigoplus\limits_{i\in \Z} \Hom_{D^b( \varLambda_Q)}(X,Y[2i]).$$

 Let $\mathrm{Proj}\varLambda_Q$ be the subcategory of all finitely generated projective $\varLambda_Q$-modules. It is well known that there is a triangle equivalence $p: D^b(\varLambda_Q)\to K^b(\mathrm{Proj}\varLambda_Q)$ taking every object to its projective resolution. Denote by $\mathrm{rep}^n\varLambda_Q$ the finite dimensional nilpotent $\varLambda_Q$-modules. Note that $\mathrm{rep}^n\varLambda_Q$ is an abelian subcategory of $\mathrm{mod}\varLambda_Q$, so we can consider the bounded derived category $D^b(\mathrm{rep}^n\varLambda_Q)$ of  $\mathrm{rep}^n\varLambda_Q$. The image $p(D^b(\mathrm{rep}^n\varLambda_Q))$ is a Hom-finite triangulated subcategory of $K^b(\mathrm{Proj}\varLambda_Q)$ since $D^b(\mathrm{rep}^n \varLambda_Q)$ is Hom-finite. 

 Let us recall Keller's construction of the triangulated hull of a root category \cite{Keller2005}.

 Denote $C_{dg}(\mathrm{Proj}\varLambda_Q)$ by the DG category of bounded complexes of finitely generated projective $\varLambda_Q$-modules. Let $\sA$ be the smallest DG subcategory of $C_{dg}(\mathrm{Proj}\varLambda_Q)$ consisting of objects of $p(D^b(\mathrm{rep}^n \varLambda_Q))$. Then
 $$D^b(\mathrm{rep}^n \varLambda_Q) \simeq  p(D^b(\mathrm{rep}^n \varLambda_Q))=H^0(\sA).$$
 Let $\sB$ be the DG category having the same objects with $\sA$, and the morphism spaces in $\sB$ are defined by
   $$\sB(X,Y):=\bigoplus\limits_{i\in \Z} \sA(X,Y[2i])$$
 Then we have that
   $$D^{b}(\mathrm{rep}^n \varLambda_Q)/[2]\simeq H^0{\sA}/[2]= H^0{\sB}.$$
 Let $\sD(\sB)$ be the derived category of DG-modules of $\sB$, Keller showed that the perfect category (for definition see Remark \ref{rem1}) $\mathrm{Per}\sB\subset \sD\sB$ is a triangulated hull of $D^b(\mathrm{rep}^n \varLambda_Q)/[2]$ such that the composition functor
  $$\pi: D^b(\mathrm{rep}^n \varLambda_Q)\to D^b(\mathrm{rep}^n \varLambda_Q)/[2]\simeq H^0(\sB) \hookrightarrow \mathrm{Per}\sB$$
 is a triangle functor.
 \begin{remark}\label{rem1}
  The perfect category $\mathrm{Per}\sB$ is defined to be the smallest full triangulated subcategory of $\sD(\sB)$ containing $H^0(\sB)$ and closed under direct factors. Note that $H^0(\sB)$ is Hom-finite, it follows that $\mathrm{Per}\sB$ is also Hom-finite.
 \end{remark}

\subsection{The relative homotopy category of 2-periodic complexes}
 Consider the relative homotopy category $K_2(\mathrm{Proj}\varLambda_Q)$ of 2-periodic complexes of finitely generated projective $\varLambda_Q$-modules. The following result comes from \cite{Peng1997}.
 \begin{theorem}
  $K_2(\mathrm{Proj}\varLambda_Q)$ is a triangulated category with the shift functor defined for complex category. And the functor $\Delta: K^b(\mathrm{Proj}\varLambda_Q)\longrightarrow K_2(\mathrm{Proj}\varLambda_Q)$ takes every complex $P^{\bullet}$ to
     $$\bigoplus_i P_{2i}\leftrightarrows \bigoplus_i P_{2i+1}$$
  is exact.
 \end{theorem}
 Since $\Delta$ commutes with the shift functor $[2]$, we have an induced functor $\Delta_2:K^b(\mathrm{Proj}\varLambda_Q)/[2]\longrightarrow K_2(\mathrm{Proj}\varLambda_Q)$. Moreover $\Delta_2$ is fully faithful (cf. \cite[Lemma 3.1]{Bridgeland2013}).

 \begin{remark}
  Following from the construction of the triangulated hull and the fact $H^0(\sB)\hookrightarrow K^b(\mathrm{Proj}\varLambda_Q)/[2]\hookrightarrow K_2(\mathrm{Proj}\varLambda_Q)$, there is a fully faithful triangle functor
    $$\iota: \mathrm{Per}\sB\longrightarrow K_2(\mathrm{Proj}\varLambda_Q).$$
 \end{remark}
  According to Remark \ref{rem1}, we will regard the triangulated hull $\mathrm{Per}\sB$ as the Hom-finite triangulated subcategory of $K_2(\mathrm{Proj}\varLambda_Q)$ in the following sections and denote it by $\sR_{\varLambda_Q}$.

  In $K_2(\mathrm{Proj}\varLambda_Q)$, we need the notion of tensor products.
 \begin{definition}
  For  $\begin{tikzcd} N^{\bullet}: N_0 \arrow[r,"d_N^0",shift left=0.3ex] &N_1\arrow[l,"d_N^1",shift left=0.3ex]\end{tikzcd}\in K_2(\mathrm{Proj}\varLambda_Q)$ and a 2-periodic complex of finite dimensional vector spaces $\begin{tikzcd} V^{\bullet}: V_0 \arrow[r,"d_V^0",shift left=0.3ex] &V_1,\arrow[l,"d_V^1",shift left=0.3ex]\end{tikzcd}$ define the tensor product $V^{\bullet}\otimes N^{\bullet}$ by
  $$\begin{tikzcd} V^{\bullet}\otimes N^{\bullet}: (V_0\otimes N_0) \oplus ( V_1\otimes N_1) \arrow[r,"d^0",shift left=0.3ex] &(V_1\otimes N_0) \oplus (V_0\otimes N_1)\arrow[l,"d^1",shift left=0.3ex]\end{tikzcd}$$
   where $d^0=\begin{bmatrix}d_V^0\otimes 1  &-1\otimes d_N^1\\ 1\otimes d_N^0  &d_V^1\otimes1 \end{bmatrix}$ and $d^1=\begin{bmatrix}d_V^1\otimes 1  &1\otimes d_N^1\\ -1\otimes d_N^0  &d_V^0\otimes1 \end{bmatrix}$.
 \end{definition}

  \begin{remark}
   It is easy to check that $V^{\bullet}\otimes N^{\bullet}$ belongs to $K_2(\mathrm{Proj}\varLambda_Q)$. In some sense, the definition of the tensor product of 2-periodic complexes is just the usual tensor product of complexes. Indeed, as $\Z_2$-graded vector spaces, $V^{\bullet}\otimes N^{\bullet}$ is the tensor product of $V_0\oplus V_1[-1]$ with $N_0\oplus N_1[-1]$, but the differential $d^0$ is given by the restriction of the differential $d_{V'}\otimes d_{N'}$, where 
     $$V':=V_0 \stackrel{d_V^0}\lrw V_1\stackrel{d_V^1}\lrw V_0,\qquad N':=N_0\stackrel{d_N^0}\lrw N_1\stackrel{d_N^1}\lrw N_0.$$
     and $d^1$ is the restriction of the differential  $d_{V''}\otimes d_{N''}$, where 
     $$V'':=V_1\stackrel{d_V^1}\lrw V_0\stackrel{d_V^0}\lrw  V_1, \qquad N'':=N_1\stackrel{d_N^1}\lrw N_0\stackrel{d_N^0}\lrw N_1.$$
  \end{remark}

  \begin{proposition} \label{prop1.6}
  For $P_{\cdot}$ and $Q_{\cdot}$ in $K^b(\mathrm{Proj}\varLambda_Q)$, we have that
   $$\Delta(P^{\bullet}\otimes Q^{\bullet})=\Delta(P^{\bullet})\otimes \Delta(Q^{\bullet}).$$
  \begin{proof}
  Without loss of generality, we may assume that $P^{\bullet}=0\to P_0\to \cdots \to P_n\to 0$ and $Q^{\bullet}=0\to Q_0\to \cdots \to Q_m\to 0$.  Then
  :
   $$\begin{tikzcd} \Delta(P^{\bullet})_{\cdot}: \bigoplus\limits_{i}P_{2i} \arrow[r,"(d_P^{2i})",shift left=0.3ex] &\bigoplus\limits_{i}P_{2i+1}\arrow[l,"(d_P^{2i+1})",shift left=0.3ex]\end{tikzcd}
    $$
 and
   $$\begin{tikzcd} \Delta(Q^{\bullet})_{\cdot}: \bigoplus\limits_{i}Q_{2i} \arrow[r,"(d^{2i}_Q)",shift left=0.3ex] &\bigoplus\limits_{i}Q_{2i+1}\arrow[l,"(d^{2i+1}_Q)",shift left=0.3ex]\end{tikzcd}.
   $$
  Then for $r\in \N$, as graded vector spaces, we have that
  $$\Delta(P^{\bullet}\otimes Q^{\bullet})_0=\bigoplus\limits_{r}\bigoplus\limits_{s+t=2r} P_s\otimes Q_t =(\bigoplus\limits_{i}P_{2i}\otimes \bigoplus\limits_j Q_{2j}) \oplus (\bigoplus\limits_{i}P_{2i+1}\otimes \bigoplus\limits_j Q_{2j+1})$$
  So $\Delta(P^{\bullet}\otimes Q^{\bullet})_0=(\Delta(P^{\bullet})\otimes \Delta(Q^{\bullet}))_0$. Similarly we can show it for $\Delta(P^{\bullet}\otimes Q^{\bullet})_1$.

  For $x_i\otimes y_j\in \Delta(P^{\bullet}\otimes Q^{\bullet})$, if $i$ and $j$ are both even, then $d_{\Delta}(x_i\otimes y_j)= d_{P^{\bullet}\otimes Q^{\bullet}}(x_i\otimes y_j)=d(x_i)\otimes y_j+(-1)^{|x_i|}x_i\otimes d(y_j)$. On the other hand,
    $$d_{\Delta\otimes \Delta}(x_i\otimes y_j)=\begin{bmatrix}d_{\Delta(P)}^0\otimes 1  &1\otimes d_{\Delta(Q)}^1\\ 1\otimes d_{\Delta(Q)}^0  &d_{\Delta(P)}^1\otimes1 \end{bmatrix} \begin{bmatrix}x_i\otimes y_j \\ 0 \end{bmatrix}=d_{\Delta(P)}^0(x_i)\otimes y_j+(-1)^{|x_i|}x_i\otimes d_{\Delta(Q)}^0(y_j).  $$
  Hence, the differentials of both sides are the same.
  \end{proof}
 \end{proposition}

\subsection{Euler form}
  Let $A$ be the path algebra $kQ$ of the acyclic quiver $Q$ and $\varLambda_Q$ the preprojective algebra of $Q$. Denote $S_i$ by the simple $A$-modules associated to $i\in Q_0$. Also, they are simple $\varLambda_Q$-modules.
  Let $<,>_A$ be the Euler form of $A$ on the Grothendieck group $K_0(A)$ of the category $\mathrm{mod} A$ of finite dimensional $A$-modules (also the derived category $D^b(\mathrm{mod} A)$). Denote by $\hat{M}$ the class of $M\in \mathrm{mod}A$ in $K_0(A)$, then
    $$<\hat{M},\hat{N}>_A:= \sum\limits_{i\in \Z} (-1)^i \dim_k \Hom_{D^b(A)}(M,N[i]).$$
  And the symmetric Euler form $(-,-)_A$ is given by
     $$(\hat{M},\hat{N})_A=<\hat{M},\hat{N}>_A+<\hat{N},\hat{M}>_A.$$
  The Euler form $<,>_{\varLambda_Q}$ of $\varLambda_Q$ is defined similarly. Note that $D^b(\mathrm{rep}^n \varLambda_Q)$ is a 2-Calabi-Yau category, we have the following result. To simplify notations, we will denote $|V|$ by the dimension of the vector space $V$.

  \begin{lemma} \label{lem1}
  For all simple modules $S_i,S_j$, $i,j\in Q_0$, we have that
          $$<\hat{S_i},\hat{S_j}>_{\varLambda_Q}=(\hat{S_i},\hat{S_j})_A.$$
  In particular, $<-,->_{\varLambda_Q}$ is a symmetric bilinear form.
  \begin{proof}
   If $i=j$, then
            $$<\hat{S_i},\hat{S_i}>_{\varLambda_Q}=|\Hom_{\varLambda_Q}(S_i,S_i)|- |\operatorname{Ext}^1_{\varLambda_Q}(S_i,S_i)| + |\operatorname{Ext}^2_{\varLambda_Q}(S_i,S_i)|=2.$$
   by direct calculations.

  If $i\neq j$, we know that $|\operatorname{Ext}^1_A(S_i,S_j)|=\#\{\alpha:i\to j| \alpha\in Q_1\}=:a_{ij}$. Note that the quiver $Q$ is acyclic, either $a_{ij}=0$ or $a_{ji}=0$. Then we have that
           $$<\hat{S_i},\hat{S_i}>_{\varLambda_Q}=|\Hom_{\varLambda_Q}(S_i,S_j)|- |\operatorname{Ext}^1_{\varLambda_Q}(S_i,S_j)| + |\operatorname{Ext}^2_{\varLambda_Q}(S_i,S_j)|=-(a_{ij}+a_{ji}).$$
 On the other hand, 
  $$|\Hom_A(S_i,S_j)|=\delta_{ij}\ and\ |\operatorname{Ext}^2_A(S_i,S_j)|=0$$
 for any $i,j$, then $<\hat{S_i},\hat{S_i}>_{\varLambda_Q}=(\hat{S_i},\hat{S_i})_A$=2. If $i\neq j$, $<\hat{S_i},\hat{S_j}>_{\varLambda_Q}=(\hat{S_i},\hat{S_j})_A= -(a_{ij}+a_{ji})$. The proof is completed.
  \end{proof}
 \end{lemma}
  Note that the composition $D^b(\mathrm{rep}^n \varLambda_Q)\stackrel{\pi}\lrw \sR_{\varLambda_Q} \subset K_2(\mathrm{Proj}\varLambda_Q)$ is a triangle functor and $K_0(D^b(\mathrm{rep}^n\varLambda_Q))\cong \Z Q_0\cong K_0(K_2(\mathrm{Proj}\varLambda_Q))$, it follows that there is an isomorphism $\eta: K_0(\sR_{\varLambda_Q})\cong \Z Q_0$ by sending $\hat{S_i}$ to $\alpha_i$.

  To construct a Kac-Moody Lie algebra of $\sR_{\varLambda_Q}$, we need a bilinear form $(-|-)_{\sR_{\varLambda_Q}}$ on the Grothendieck group $K_0(\sR_{\varLambda_Q})$.  Define a bilinear form on $K_0(\sR_{\varLambda_Q})$ by
  $$(\hat{M}|\hat{N})_{\sR_{\varLambda_Q}}:=\dim_k \Hom_{\sR_{\varLambda_Q}}(M,N)-\dim_k \Hom_{\sR_{\varLambda_Q}}(M,N[1]).$$
  for $\hat{M},\hat{N}$ in $K_0(\sR_{\varLambda_Q})$.

  \begin{remark}\label{rmk1}
   Note that $\sR_{\varLambda_Q}$ is 2-periodic, the bilinear form above is well-defined, i.e. for any triangle $N'\to N\to N''\to N'[1]$ in $\sR_{\varLambda_Q}$ and $M\in \sR_{\varLambda_Q}$, we have
          $$(\hat{M}|\hat{N})_{\sR_{\varLambda_Q}}:=(\hat{M}|\hat{N'})_{\sR_{\varLambda_Q}}+(\hat{M}|\hat{N''})_{\sR_{\varLambda_Q}}.$$
   Moreover, the bilinear form is symmetric. Indeed, by using the fact 
     $$\Hom_{\sR_{\varLambda_Q}}(S_i,S_j)=\Hom_{D^b(\varLambda_Q)}(S_i,S_j)\oplus \Hom_{D^b(\varLambda_Q)}(S_i,S_j[2]),$$ 
 and 
     $$\Hom_{\sR_{\varLambda_Q}}(S_i,S_j[1])=\Hom_{D^b(\varLambda_Q)}(S_j,S_i[1]),$$
  it follows that
            $$(\hat{S_i}|\hat{S_j})_{\sR_{\varLambda_Q}}= <\hat{S_i},\hat{S_j}>_{\varLambda_Q}$$
    and $<-,->_{\varLambda_Q}$ is symmetric by Lemma \ref{lem1}.
  \end{remark}
\begin{corollary}
   The triangle functor $\pi: D^b(\mathrm{rep}^n\varLambda_Q)\lrw \sR_{\varLambda_Q}$ induces an isotropy $\pi^*: K_0(\varLambda_Q)\lrw K_0(\sR_{\varLambda_Q})$.
 \begin{proof}
   We know that the map $\pi^*(\hat{S_i})=\hat{S_i}$ is an isomorphism of groups by Remark \ref{rmk1}. Also Remark \ref{rmk1} tells us that $<\hat{S_i},\hat{S_j}>_{\varLambda_Q}=(\hat{S_i}|\hat{S_j})_{\sR_{\varLambda_Q}}$.
 \end{proof}
\end{corollary}

\subsection{Integral Hall algebra}
   Fix a finite field $k=\F_q$. Let $R$ be a finitary algebra over $k$ and $\mathrm{mod}R$ the category of finite dimensional right $R$-modules. The integral Ringel–Hall algebra $\sH(R)$ of $R$ is by definition the free abelian group with the basis consisting of all isoclasses $[M]$ for $M\in \mathrm{mod}R$. The multiplication is given by
      $$[M][N]:=\sum_{[L]}g_{MN}^L [L],$$
   where $g_{MN}^L$ is the number of submodules $X$ of $L$ such that $X\cong N$ and $L/X\cong M$.

  Denote by $C\sH(R)$ the subalgebra of $\sH(R)$ generated by $[S_i]$, $i\in Q_0$, which is the so-called composition algebra of $\sH(R)$. The following proposition is given in Proposition $5.1$ of \cite{Chen2008}.
 \begin{proposition}\label{prop1}
    Let $R$ be a finitary $k$-algebra and let $R'$ be a factor algebra of $R$. Then there
   are epimorphisms $\sH(R)\to \sH(R')$ and $C\sH(R)\to C\sH(R')$ of $\Z$-algebras.
 \end{proposition}

\subsection{The Ringel-Hall Lie algebra}
  Recall the definition of the Ringel-Hall Lie algebra of a root category $\sR$ following \cite{Peng2000}. By $\mathrm{ind}\sR$ we denote the set of all isoclasses of indecomposable objects in $\sR$.

   Given $X$, $Y, L\in \mathrm{Obj}(\sR)$, we define
   \begin{center}
       $W(X,Y;L):=\{(f,g,h)\in \Hom_{\sR}(Y,L)\times \Hom_{\sR}(L,X)\times \Hom_{\sR}(X,Y[1])|$ $Y\stackrel{f}\to L\stackrel{g}\to X \stackrel{h}\to Y[1]$ is a triangle $\}$.\\
       $V(X,Y;L):=W(X,Y;L)/(\mathrm{Aut}(Y)\times \mathrm{Aut}(X)).$
       \end{center}
   where the action of $\mathrm{Aut}(Y)\times \mathrm{Aut}(X)$ on $W(X,Y;L)$ is given by
    $$(a,c).(f,g,h)=(fa,c^{-1}g, a[1]^{-1}hc)$$
    for $(f,g,h)\in W(X,Y;L)$ and $(a,c)\in \mathrm{Aut}(X)\times \mathrm{Aut}(Y)$. In the following, we will denote  by $(f,g,h)^{\hat{}}$ the orbit of $(f,g,h)$ in $V(X,Y;L)$.

    Let $\Hom_{\sR}(Y,L)_X$ be the subset of  $\Hom_{\sR}(Y,L)$ consisting of morphisms $f:Y\to L$ such that $\mathrm{cone}f\cong X$. Consider the action of $\mathrm{Aut}(Y)$ on $\Hom_{\sR}(Y,L)_X$ by $a.f=fa$, the orbit space is denoted by
       $$\Hom_{\sR}(Y,L)^*_X:=\Hom_{\sR}(Y,L)_X/\mathrm{Aut}(Y).$$
   Dually, consider the subset $\Hom_{\sR}(L,X)_{Y[1]}$ of $\Hom_{\sR}(L,X)$ and $\mathrm{Aut}(X)$ acts on $\Hom_{\sR}(L,X)_{Y[1]}$ by $c.g=gc^{-1}$,  we have another orbit space $\Hom_{\sR}(L,X)^*_{Y[1]}$.

  \begin{proposition}[\cite{Xiao2008}]
    $ |V(X,Y;L)|=|\Hom_{\sR}(Y,L)^*_X|=|\Hom_{\sR}(L,X)^*_{Y[1]}|$.
  \end{proposition}

   Let $\ssh$ be the subgroup of $K_0(\sR)\otimes_{\Z} \Q$ generated by $\frac{h_{M}}{d(M)}$ with $M\in \mathrm{ind}\sR$ and $d(M)=\dim_k \operatorname{End}M/\mathrm{rad}(\operatorname{End}M)$. One can naturally extend the symmetric Euler form $(-|-)_{\sR}$ to $\ssh\times \ssh$. Let $\ssn$ be the free abelian group with the basis $\{u_X|X\in \mathrm{ind}\sR\}$. Let
       $$\ssg(\sR)=\ssh\oplus \ssn$$
   be a direct sum of $\Z$-modules. The Lie operation is given as follows:
   \begin{itemize}
       \item for any indecomposable objects $X,Y\in \sR$,
          $$[u_X, u_Y]=\sum_{L\in \mathrm{ind}\sR} (F_{XY}^L-F_{YX}^L)u_L-\delta_{X,Y[1]}\frac{h_X}{d(X)}. $$
       where $F_{XY}^L:=|V(X,Y;L)|$.
       \item for any objects $X,Y\in \sR$ with $Y$ indecomposable
           $$[h_{X},u_Y]=(h_X|h_Y)_{\sR}u_Y = -[u_Y,h_X]$$
       \item $[\mathfrak{h},\mathfrak{h}]=0$.
   \end{itemize}
   Consider the quotient algebra
       $$\ssg(\sR)_{(q-1)}= \ssg(\sR)/(q-1)\ssg(\sR).$$
   Then by Peng-Xiao \cite{Peng2000} we know that $\ssg(\sR)_{(q-1)}$ is a Lie algebra over $\Z/(q-1)$.

   Recall that $A$ is the path algebra $kQ$ of the acyclic quiver $Q$ of non-Dynkin type and $\varLambda_Q$ is the preprojective algebra. Since $A$ is hereditary, the root category $\sR_A$ is triangle equivalent to the homotopy category $K_2(\mathrm{Proj}A)$ of 2-periodic complexes of finitely generated projective $A$-modules.  Recall that $\sR_{\varLambda_Q}$ is the triangulated hull of $D^b(\mathrm{rep}^n\varLambda_Q)/[2]$. Since $\sR_A$ and $\sR_{\varLambda_Q}$ are both Hom-finite  $k$-linear triangulated category with the shift functor $[1]$, we have the corresponding Ringel-Hall Lie algebras $\ssg(\sR_A)_{q-1}=\ssh(A) \oplus \ssn(A)$ and $\ssg(\sR_{\varLambda_Q})_{q-1}$ over $\Z/(q-1)$.

   Set $h_i:=\frac{h_{S_i}}{d(S_i)}$, consider the Lie subalgebra $\ssg_0(\sR_A)$ (resp. $\ssg_0(\sR_{\varLambda_Q})$) of $\ssg(\sR_A)_{(q-1)}$ (resp. $\ssg(\sR_{\varLambda_Q})_{(q-1)}$) generated by $\{u_{S_j}, u_{S_j[1]}, h_j| j\in Q_0\}$ (resp. $\{\hat{u}_{S_j}, \hat{u}_{S_j[1]},h_j| j\in Q_0\}$), and $\ssn^+_0(A)$ (resp. $\ssn^+_0(\sR_{\varLambda_Q})$) generated by $\{u_{S_j}|j\in Q_0\}$ (resp. $\{\hat{u}_{S_j}|j\in Q_0\}$). By $\hat{u}_{S_j}$ we mean $\hat{u}_{\Delta\cdot p(S_j)}$, i.e., the image of $S_j$ in $H^0(\sB)\subset K_2^b(\mathrm{Proj}\varLambda_Q))$.

\section{Spherical Twist Functors and Weyl Groups}
 \subsection{Construction of spherical twist functors}
       Since $D^b(\mathrm{rep}^n\varLambda_Q)$ is a 2-Calabi-Yau category, every simple object is a 2-spherical object (i.e., $\mathbb{R}\Hom(S_i,S_i)\cong k\oplus k[-2]$.) By \cite{Seidel2001}, for a simple $\varLambda_Q$-module $S$,  let $P_S$ be the minimal  projective resolution of $S$, then the spherical twist functor $T_S: D^b(\mathrm{rep}^n\varLambda_Q)\to D^b(\mathrm{rep}^n\varLambda_Q)$ is given by
           $$T_S(X)=cone(\Hom^{\bullet}(P_S,X)\otimes P_S \stackrel{ev}\longrightarrow X)$$
       where $\Hom^{\bullet}(P_S,X)$ is the usual complex of vector spaces and $ev$ is the natural evaluation map, which is a cochain map. For any chain map $f:X\to Y$,
            $$T_S(f)=\begin{bmatrix}f &0\\0 &f_{*}\otimes 1 \end{bmatrix}.$$
       Moreover, $T_S$ is a triangle autoequivalence on $D^b(\mathrm{rep}^n\varLambda_Q)$ (cf. \cite[Proposition 2.10]{Seidel2001}).
       \begin{remark}\ \
   
       (i) The spherical twist  functor $T_S$ is well defined. Note that $\Hom^{\bullet}(P_S,X)\otimes P_S$ is a bounded complex since both $P_S$ and $X$ are bounded. Furthermore, we have that $\operatorname{H}^*(\Hom^{\bullet}(P_S,X)\otimes P_S)\cong \operatorname{H}^*(\Hom^{\bullet}(P_S,X))\otimes \operatorname{H}^*(P_S)$ has finite total dimension, so $\Hom^{\bullet}(P_S,X)\otimes P_S$ belongs to $D^b(\mathrm{rep}^n\varLambda_Q)$, it follows that the cone $T_S(X)$ lies in $D^b(\mathrm{rep}^n\varLambda_Q)$.
   
       (ii)The spherical twist functor $T_S:D^b(\mathrm{rep}^n\varLambda_Q)\to D^b(\mathrm{rep}^n\varLambda_Q)$ is induced by a DG functor of some DG categories. Indeed, we have known that $\sA$ is the DG enhancement of $D^b(\mathrm{rep}^n\varLambda_Q)$, that is $\operatorname{H}^0(\sA)\simeq D^b(\mathrm{rep}^n\varLambda_Q)$. Moreover, there is a DG functor
           $$T'_S: \sA\lrw \sA,\ P_X\mapsto cone(ev:\Hom^{\bullet}(P_S,P_X)\otimes P_S\to P_X),$$
         sending $f\in \Hom^i(P_X,P_Y)$ to
         $$T'_S(f):= \begin{bmatrix}f &0\\0 &(-1)^if_{*}\otimes 1 \end{bmatrix}.$$
        Therefore, if we regard $D^b(\mathrm{rep}^n\varLambda_Q)$ as $H^0(\sA)$, then $T_S=H^0(T'_S)$.
       \end{remark}
   
       \begin{example}\label{ex3.2}
        For $i\in Q_0$, let $S_i$ be the simple $\varLambda_Q$-module. The minimal projective resolution $P_{S_i}$ of $S_i$ is given by
        $$P_{S_i}:=0\lrw P_i\lrw \bigoplus\limits_{t(h)=i}P_{s(h)}\lrw P_i\lrw 0$$
        Then $\Hom^{\bullet}(P_{S_i},S_i)=k\oplus k[-2]$, and $ev=(f_0,f_2):P_{S_i}\oplus P_{S_i}[-2]\to S_i$, where $f_0\in \Hom_{\varLambda_Q}(S_i,S_i)$ and $f_1\in \operatorname{Ext}^2_{\varLambda_Q}(S_i,S_i)$ are nontrivial. Hence $T_{S_i}(S_i)=cone(ev)\cong S_i[-1]$.
   
        Set $a_{ji}=\dim\operatorname{Ext}^1_{\varLambda_Q}(S_i,S_j)$. Since $\Hom^{\bullet}(P_{S_i},S_j)\otimes P_{S_i}=P_{S_i}^{\oplus a_{ji}}[-1]$, it follows that $T_{S_i}(S_j)=I_{ji}$, where $I_{ji}$ is an extension of $S_i^{\oplus a_{ji}}$ by $S_j$, i.e., $I_{ji}=(V_i\oplus V_j, h\in \overline{Q}_1)$, where $V_i=k^{\oplus a_{ji}}$, $V_j=k$, $h_l=(0,\cdots,0,1,0,\cdots,0):V_i\to V_j$ with 1 in the $l$-th term,  if $t(h_l)=j$, $1\leq l\leq a_{ji}$, and $\bar{h}_l=0$ for all $l$. Note that $I_{ji}$ is indecomposable.
       \end{example}
   
       We recall some useful properties for spherical twist functors from \cite[Proposition 6.8]{Shiraishi2016a}.
       \bproposition
       Let $S$ be a spherical object in $D^b(\mathrm{rep}^n\varLambda_Q)$.
   
       (i)For an integer $l\in \Z$, we have that $T_{S[l]}\cong T_S$.
   
       (ii) We have that $T_S(S)=S[-1]$.
   
       \eproposition
   
        Note that $T'_S$ commutes with shift functor $[2]$, we have an induced functor $T'_{2,S}: \sB\to \sB$.  We have the following diagram commute.
        $$\begin{tikzcd} 
         D^b(\mathrm{rep}^n\varLambda_Q)\arrow[r,"T_S"]\arrow[d,"p"]  & D^b(\mathrm{rep}^n\varLambda_Q)\arrow[d,"p"] \\
         H^0(\sA)\arrow[r,"H^0(T'_S)"]\arrow[d,"F"]       &H^0(\sA)\arrow[d,"F"] \\
         H^0(\sB)\arrow[r,"H^0(T'_{2,S})"]    &H^0(\sB). \end{tikzcd}\\$$
       Denote by $T_{2,S}$ the induced functor $H^0(T'_{2,S})$, then $T_{2,S} : H^0(\sB)\to H^0(\sB)$ is an equivalence.
   
     To construct a functor $\bar{T}_{2,S}: K_2(\mathrm{Proj}\varLambda_Q) \to K_2(\mathrm{Proj}\varLambda_Q)$ extending the induced twist functor $T_{2,S}: H^0(\sB)\to H^0(\sB)$, for $P^{\bullet}$, $Q^{\bullet}\in K_2(\mathrm{Proj}\varLambda_Q))$, we define $B_2(P^{\bullet},Q^{\bullet})$ to be a 2-periodic complexes of vector spaces  as follows:
     $$\begin{tikzcd} B_2(P^{\bullet},Q^{\bullet}):  Hom_{\Lambda_Q}(P^{\bullet},Q^{\bullet}) \arrow[r,"d^0",shift left=0.3ex] & Hom_{\Lambda_Q}(P^{\bullet},Q^{\bullet}[1])\arrow[l,"d ^1",shift left=0.3ex]\end{tikzcd}$$
     where $\Hom_{\varLambda_Q}(P^{\bullet},Q^{\bullet})=\Hom_{\varLambda_Q}(P_0,Q_0)\oplus \Hom_{\varLambda_Q}(P_1,Q_1)$,  $d^0$ is defined by
         $$d^0(f_0)=d^0_Q(f_0)-f_0d^1_P,\qquad d^0(f_1)=d^1_Q(f_1)-f_1d^0_P,$$
      for $f=(f_0,f_1)\in \Hom_{\varLambda_Q}(P^{\bullet},Q^{\bullet})$,  and $d^1$ is given by
         $$d^1(g_0)=d^1_Q(g_0)+g_0d^1_P,\qquad d^1(g_1)=d^0_Q(g_1)+g_1d^0_P,$$
     for $g=(g_0,g_1)\in \Hom_{\varLambda_Q}(P^{\bullet},Q^{\bullet}[1])$.  It can be checked that $d^0d^1=0$ and $d^1d^0=0$.  Note if $P^{\bullet}$ and $Q^{\bullet}$ are in the image of the compression functor $\Delta: K^b(\mathrm{Proj}\varLambda_Q)\to K_2(\mathrm{Proj}\varLambda_Q)$, then $B_2(P^{\bullet},Q^{\bullet})$ is precisely the compression of $\Hom^{\bullet}(P^{\bullet},Q^{\bullet})$ by definition.
   
     Define the functor $\bar{T}_{2,S}: K_2(\mathrm{Proj}\varLambda_Q)\to K_2(\mathrm{Proj}\varLambda_Q)$ by
        $$\bar{T}_{2,S} (Q^{\bullet})= Cone(ev: B_2(\Delta(P_S),Q^{\bullet})\otimes \Delta(P_S)\longrightarrow Q^{\bullet}).$$
     for any object $Q^{\bullet}\in K_2(\mathrm{Proj}\varLambda_Q)$, 
        $$\bar{T}_{2,S}(f):=\begin{bmatrix} f &0\\0 &f_*\otimes 1\end{bmatrix}:cone(ev)\lrw cone(ev'),$$
     for any $f\in \Hom_{K^2(\mathrm{Proj}\varLambda_Q)}(Q^{\bullet},Q'^{\bullet})$. Because $f_{*}\otimes 1: B_2(\Delta(P_S),Q^{\bullet})\otimes \Delta(P_S) \to  B_2(\Delta(P_S),Q'^{\bullet})\otimes \Delta(P_S)$ is also a morphism of 2-periodic complexes and $f\circ ev=ev'\circ (f_*\otimes 1)$,  $\bar{T}_{2,S}(f):cone(ev)\to cone(ev')$ is well defined.
   
     \blemma\label{lem6}
     The functor $\bar{T}_{2,S}:K_2(\mathrm{Proj}\varLambda_Q)\to K_2(\mathrm{Proj}\varLambda_Q)$ restricts to a functor $T_{2,S}:H^0(\sB)\to H^0(\sB)$ along the fully faithful functor $\Delta_2: K^b(\mathrm{Proj}\varLambda_Q)/[2]\lrw K_2(\mathrm{Proj}\varLambda_Q)$.
     \bproof
      By the definitions of $\bar{T}_{2,S}$ and $T_{2,S}$, they are the same in morphism spaces. It suffices to show $\bar{T}_{2,S}(\Delta_2(P^{\bullet}))=\Delta_2(T_{2,S}(P^{\bullet}))=\Delta(T_{2,S}(P^{\bullet}))$ for each object $P^{\bullet}$ in $H^0(\sB)$. Note that $\Delta(\Hom^{\bullet}(P_S,P^{\bullet})\otimes P_S)=\Delta(\Hom^{\bullet}(P_S,P^{\bullet}))\otimes \Delta(P_S)=B_2(\Delta(P_S),\Delta(P^{\bullet}))\otimes \Delta(P_S)$ by Proposition \ref{prop1.6}, then the evaluation map 
        $$ev':B_2(\Delta(P_S),\Delta(P^{\bullet}))\otimes \Delta(P_S)\to \Delta(P^{\bullet})$$
      in $K_2(\mathrm{Proj}\varLambda_Q)$ is precisely $\Delta(ev:\Hom^{\bullet}(P_S,P^{\bullet})\otimes P_S\to P^{\bullet})$.  Thus we have $\bar{T}_{2,S}(\Delta(P^{\bullet}))=cone(ev')=cone(\Delta(ev))=\Delta(cone(ev))=T_{2,S}(P^{\bullet})$, which completes the proof.
     \eproof
     \elemma
   
     Next, we will prove the following proposition. Our strategy is firstly to show $\bar{T}_{2,S}:K_2(\mathrm{Proj}\varLambda_Q)\lrw K_2(\mathrm{Proj}\varLambda_Q)$ is a triangle functor. Due to the above Lemma \ref{lem6} and the fact that $\sR_{\varLambda_Q}$ is the smallest triangulated subcategory of $K_2(\mathrm{Proj}\varLambda_Q)$ containing $H^0(\sB)$, it follows that $\bar{T}_{2,S}$ sends $\sR_{\varLambda_Q}$ to $\sR_{\varLambda_Q}$.
     \bproposition\label{prop2.3}
     The restriction functor $\bar{T}_{2,S}|_{\sR_{\varLambda_Q}}: \sR_{\varLambda_Q} \to \sR_{\varLambda_Q}$ of the functor $\bar{T}_{2,S}$ along the triangulated hull $\sR_{\varLambda_Q}$ is a triangle equivalence.
     \eproposition
     In order to prove Proposition \ref{prop2.3}, we introduce a DG category $\sC$:
       \begin{center}$Ob(\sC)$=\text{all 2-periodic complexes of  finitely generated projective} $\varLambda_Q$\text{-modules}\end{center}
      For $M^{\bullet},N^{\bullet}\in \sC$, the morphism space is defined by
       $$\Hom_{\sC}^{\bullet}(M^{\bullet},N^{\bullet}):=\bigoplus\limits_{i\in \Z_2} \Hom_{\varLambda_Q}(M^{\bullet},N^{\bullet}[i]),$$
     where $\Hom_{\varLambda_Q}(M^{\bullet},N^{\bullet})=\Hom_{\varLambda_Q}(M_0,N_0)\oplus \Hom_{\varLambda_Q}(M_1,N_1)$, and the differential $d^{\cdot}$ of $\Hom_{\sC}^{\bullet}(M^{\bullet},N^{\bullet})$ is given by
       $$d(f):=d_N f-(-1)^{i} fd_M,$$
     for $f\in \Hom^i_{\sC}(M^{\bullet},N^{\bullet})$. The composition map $\Hom_{\sC}^{\bullet}(N^{\bullet},L^{\bullet})\otimes_k \Hom_{\sC}^{\bullet}(M^{\bullet},N^{\bullet})\to \Hom_{\sC}^{\bullet}(M^{\bullet},L^{\bullet})$ is given by $f^i\otimes g^j\mapsto f^ig^j$ for $i,j=\bar{0},\bar{1}$. It can be checked that $\sC$ is a DG category. And it is clear thar $H^0(\sC)=K_2(\mathrm{Proj}\varLambda_Q)$.
   
     Define a functor $\bar{T_2}: \sC \to \sC$ by
        $$\bar{T_2}(M^{\bullet}):=cone(ev:B_2(\Delta(P_S),M^{\bullet})\otimes \Delta(P_S)\to M^{\bullet}),$$
     for $f\in \Hom^i_{\sC}(M^{\bullet},N^{\bullet})$,
        $$\bar{T_2}(f):= \begin{bmatrix}f &0\\0 &(-1)^{i}f_{*}\otimes 1 \end{bmatrix}\in \Hom_\sC(\bar{T_2}(X),\bar{T_2}(Y))^i.$$
   
     \blemma
      The functor $\bar{T_2}:\sC\to \sC$ defined as above is a DG functor. In particular, the homotopy functor $H^0(\bar{T}_{2})$ is the functor $\bar{T}_{2,S}: K_2(\mathrm{Proj}\varLambda_Q)\to K_2(\mathrm{Proj}\varLambda_Q)$.
     \bproof
      We need to show $(\bar{T_2})_{X,Y}: \Hom_{\sC}^{\bullet}(X,Y)\to \Hom_{\sC}^{\bullet}(\bar{T_2}(X),\bar{T_2}(Y))$ is a strict morphism for any $X,Y\in \mathrm{Obj}(\sC)$, i.e., $d_{\Hom^{\bullet}(\bar{T_2}(X),\bar{T_2}(Y))}(\bar{T_2})_{X,Y}=(\bar{T_2})_{X,Y}d_{\Hom^{\bullet}(X,Y)}$.
      Namely, for any $f\in \Hom^i_{\sC}(X,Y)$,
        $$d_{\bar{T_2}(X)}\bar{T_2}(f)-(-1)^{i}\bar{T_2}(f)d_{\bar{T_2}(X)}=\begin{bmatrix}d(f) &0\\0 &(-1)^{i+1}(df)_{*}\otimes 1 \end{bmatrix}$$
      where $d(f)=d_{Y}f-(-1)^i fd_{X}.$
       $$\begin{aligned}
        &\mathrm{LHS}\\
        &=\begin{bmatrix}d_Y &ev\\0 &-d_{B_2(\Delta(P_S),Y)\otimes \Delta(P_S)} \end{bmatrix}\begin{bmatrix}f &0\\0 &(-1)^{i}f_{*}\otimes 1 \end{bmatrix}\\
        &-(-1)^{i}\begin{bmatrix}f &0\\0 &(-1)^{i}f_{*}\otimes 1 \end{bmatrix}\begin{bmatrix}d_X &ev\\0 &-d_{B_2(\Delta(P_S),X)\otimes \Delta(P_S)} \end{bmatrix} \\
        &= \begin{bmatrix}d_Yf-(-1)^i fd_X &(-1)^i [ev\cdot(f_*\otimes 1)-f\cdot ev]\\0 &(f_{*}\otimes 1)d_{B_2(\Delta(P_S),X)\otimes \Delta(P_S)}-(-1)^id_{B_2(\Delta(P_S),Y)\otimes \Delta(P_S)}(f_*\otimes 1) \end{bmatrix} \\
        &=\begin{bmatrix}d_Yf-(-1)^i fd_X &0\\0 &(fd_X)_*\otimes 1-(-1)^i (d_Yf)_*\otimes 1 \end{bmatrix}\\
        &=\mathrm{RHS}.
       \end{aligned}$$
      The second statement is clear from the definition of $\bar{T}_2$.
     \eproof
     \elemma
      So $\bar{T}_{2,S}:K_2(\mathrm{Proj}\varLambda_Q) \to K_2(\mathrm{Proj}\varLambda_Q)$ is a triangle functor by \cite[Theorem 4.4.41]{Yang2019}, which implies the restriction $\bar{T}_{2,S}|_{\sR_{\varLambda_Q}}: \sR_{\varLambda_Q}\to \sR_{\varLambda_Q}$ is also a triangle functor.
   
      \ \ \\
     \textbf{Proof of Proposition \ref{prop2.3}}
     \bproof
      We only need to show $\bar{T}_{2,S}|_{\sR_{\varLambda_Q}}: \sR_{\varLambda_Q}\to \sR_{\varLambda_Q}$ is an equivalence. Note that $\bar{T}_{2,S}|_{H^0{\sB}}=T_{2.S}$ and $T_{2,S}$ is an equivalence, hence $\bar{T}_{2,S}|_{\sR_{\varLambda_Q}}$ is fully faithful and dense on $H^0\sB$. By the definition of the triangulated hull $\sR_{\varLambda_Q}=\mathrm{Per}\sB$ (which is the smallest triangulated subcategory containing $H^0\sB$ and closed under direct summands), it follows that the functor $\bar{T}_{2,S}|_{\sR_{\varLambda_Q}}$ is fully faithful by using Five lemma continuously. It can be showed in the same way that $\bar{T}_{2,S}|_{\sR_{\varLambda_Q}}$ is dense.
     \eproof
     In the following sections, we will simply denote by $\bar{T}_{2,S}$ the restriction functor $\bar{T}_{2,S}|_{\sR_{\varLambda_Q}}$, and call it the spherical twist functor.
   
     \subsection{Isomorphisms of Lie algebras induced by spherical twist  functors}
        Fix a vertex $i\in Q_0$,  the spherical twist functor $\bar{T}_{2,S_i}:\sR_{\varLambda_Q}\to \sR_{\varLambda_Q}$ (for simplicity, denote $\bar{T}_i:=\bar{T}_{2,S_i}$) naturally induces a morphism
             $$\Phi_i: \ssg(\sR_{\varLambda_Q})_{(q-1)}\to \ssg(\sR_{\varLambda_Q})_{(q-1)},\ u_X\mapsto u_{\bar{T}_i(X)},\ h_{Y}\mapsto h_{\bar{T}_i(Y)}.$$
   
        \bproposition
          For any simple $\varLambda_Q$-module $S_i$, the spherical twist functor $\bar{T}_{2,S_i}:\sR_{\varLambda_Q}\to \sR_{\varLambda_Q}$ induces an isomorphism $\Phi_i: \ssg(\sR_{\varLambda_Q})_{q-1}\to \ssg(\sR_{\varLambda_Q})_{(q-1)}$ of Lie Algebras.
        \bproof
          First, we need to show that $\Phi_i$ is a homomorphism of Lie algebras.
   
          For $X,Y\in \mathrm{ind}\sR_{\varLambda_Q}$,  $\Phi_i([\hat{u}_X,\hat{u}_Y])=\sum_{L} (F_{XY}^L-F_{YX}^L)\hat{u}_{\bar{T}_iL}+\delta_{X,Y[1]}\frac{h_{\bar{T}_iX}}{d(X)}$. On the other hand, $[\hat{u}_{\bar{T}_iX},\hat{u}_{\bar{T}_iY}]=\sum_{L'} (F_{\bar{T_i}X\bar{T_i}Y}^{L'}-F_{\bar{T_i}Y\bar{T_i}X}^{L'})\hat{u}_{L'}+ \delta_{\bar{T}_iX,\bar{T}_iY[1]}\frac{h_{\bar{T}_iX}}{d(\bar{T_i}X)}$. Since $\bar{T}_i:\sR_{\varLambda_Q} \to \sR_{\varLambda_Q}$ is a triangle equivalence,  we have $Y\stackrel{f}\to L\stackrel{g}\to X\stackrel{h}\to Y[1]$ is a triangle iff $\bar{T}_iY\stackrel{\bar{T}_i(f)}\to \bar{T}_iL\stackrel{\bar{T}_i(g)}\to \bar{T}_iX\stackrel{\bar{T}_i(h)}\to \bar{T}_iY[1]$ is a triangle. Moreover,
          \begin{center}
            $(f,g,h)\stackrel{\mathrm{Aut}(Y)}\sim (f',g',h')$ iff \\
             $(\bar{T}_i(f),\bar{T}_i(g),\bar{T}_i(h))\stackrel{\mathrm{Aut}(\bar{T}_iY)}\sim (\bar{T}_i(f'),\bar{T}_i(g'),\bar{T}_i(h'))$,
         \end{center}
          for $\bar{T}_i$ preserves isomorphisms and $\bar{T}_i(id)=id$. Therefore, we have
                   $$|\frac{\Hom_{\sR_{\varLambda_Q}}(Y,L)_X}{\mathrm{Aut}(Y)}|=|\frac{\Hom_{\sR_{\varLambda_Q}}(\bar{T}_iY,\bar{T}_iL)_{\bar{T}_iX}}{\mathrm{Aut}(\bar{T}_iY)}|,$$
          which implies $F_{XY}^L=F_{\bar{T}_i,X\bar{T}_iY}^{\bar{T}_iL}$, it follows that
                 $$[\hat{u}_{\bar{T}_iX},\hat{u}_{\bar{T}_iY}]=\sum_{L'\cong\bar{T}_iL} (F_{XY}^{L'}-F_{YX}^{L'})\hat{u}_{\bar{T}_iL}+ \delta_{\bar{T}_iX,\bar{T}_iY[1]}\frac{h_{\bar{T}_iX}}{d(\bar{T_i}X)}.$$
          Note $X\cong Y$ iff $\bar{T}_iX\cong \bar{T}_iY$, then $\Phi_i[\hat{u}_X,\hat{u}_Y]=[\hat{u}_{\bar{T}_iX},\hat{u}_{\bar{T}_iY}]=[\Phi_i(\hat{u}_X),\Phi_i(\hat{u}_Y)]$.
   
          For $X, Y\in \sR_{\varLambda_Q}$, $Y$ is indecomposable, $\Phi([h_X,\hat{u}_Y])=(h_X|h_Y)_{\sR_{\varLambda_Q}}\hat{u}_{\bar{T}_iY}$. On the other hand, $[h_{\bar{T}_iX},\hat{u}_{\bar{T}_iY}]=(h_{\bar{T}_iX}|h_{\bar{T}_iY})_{\sR_{\varLambda_Q}}\hat{u}_{\bar{T}_iY}$. Since the triangle equivalence $\bar{T}_i$ induces an isotropy on the Grothendieck group $K_0(\sR_{\varLambda_Q})$ with respect to $(-|-)_{\sR_{\varLambda_Q}}$, we have $(h_X|h_Y)_{\sR_{\varLambda_Q}}=(h_{\bar{T}_iX}|h_{\bar{T}_iY})_{\sR_{\varLambda_Q}}$, which gives $\Phi([h_X,\hat{u}_Y])=[h_{\bar{T}_iX},\hat{u}_{\bar{T}_iY}]$.
   
          So $\Phi_i$ is a homomorphism of Lie algebras. Then it is an isomorphism and the inverse is induced by the quasi-inverse of $\bar{T}_i$.
        \eproof
        \eproposition
   
        Recall that $\ssg_0(\sR_{\varLambda_Q})$ is the Lie subalgebra of $\ssg(\sR_{\varLambda_Q})_{(q-1)}$ generated by $\hat{u}_{S_j}$, $\hat{u}_{S_j[1]}$ and $h_j$ for all $j\in Q_0$. From Example \ref{ex3.2}, we know that $T_{2,S_i}(S_i)=S_i[-1]$ and $T_{2,S_i}(S_j)=I_{ji}$ if $j\neq i$ for simple $\varLambda_Q$-module $S_j$. Then $\Phi_i:\ssg_0(\sR_{\varLambda_Q})\to \ssg(\sR_{\varLambda_Q})$ is given by
           $$\begin{aligned}
           \Phi_i(\hat{u}_{S_i})&= \hat{u}_{S_i[1]}, \qquad   \Phi_i(\hat{u}_{S_i[1]})= \hat{u}_{S_i},\\
           \Phi_i(\hat{u}_{S_j})&= u_{I_{ji}},\ \ \qquad \Phi_i(\hat{u}_{S_j[1]})= \hat{u}_{I_{ji}[1]},\ \ if\  j\neq i.\\
           \end{aligned}$$
        When $i\in Q_0$ is a source, \cite[Theorem 2.1]{Xiao2005} stated that there exists a functor $R(S_i^-)$ on root categories inducing an isomorphism $\tilde{\phi_i}$ of Kac-Moody Lie algebras, which is a lifting of Weyl group actions on root system. Namely, $\tilde{\phi}(e_i)=-f_i$, and for $j\neq i$, $\tilde{\phi_i}(e_j)=ad(e_i)^{a_{ij}}(e_j)$, etc. However, as for $\Phi_i$, this does not hold.
        \bremark
        Even if $i\in Q_0$ is a source, $\Phi_i$ still can not be restricted to an automorphism of the Lie subalgebra $\ssg_0(\sR_{\varLambda_Q})$. Indeed,
              $$[I_{ji}]\neq  \sum_{r=0}^{a_{ij}}(-1)^{r}[S_i]^{(a_{ji}-r)}[S_j][S_i]^{(r)}=[I_{ji}]-[P_{ji}].$$
        where $P_{ji}\in \mathrm{rep}^n\varLambda_Q$ is an extension of $S_i^{\oplus a_{ji}}$ by $S_j$, i.e. $P_{ji}=(V_j\oplus V_i, h\in \overline{Q}_1)$  with $V_j=k$, $V_i=k^{a_{ji}}$, $h_l=(0,\cdots,0,1,0,\cdots,0)^t:V_j\to V_i$ if $s(h_l)=j$, $1\leq l\leq a_{ji}$. Thus, $\Phi_i(S_j)\not\in \ssg_0(\sR_{\varLambda_Q})$.
        \eremark
   
        \bremark\label{rem 3.9}
         The spherical twist functor $T_{S_i}$ on $D^b(\mathrm{rep}^n\varLambda_Q)$ is the right derived functor of the reflection functor $\Sigma_i$ defined in \cite[Section 2.2]{Baumann2012a}. Indeed, it suffices to show that $H^0(T_{S_i}(M))=\Sigma_i(M)$ for any $M\in \mathrm{rep}^n\varLambda_Q$. Denoted by $\mathrm{hd}_i M$ the $S_i$-isotropic component of the head of $M\in \mathrm{rep}^n\varLambda_Q$. For $M=(\bigoplus_{j\in Q_0}M_j, M_h)$ such that $\mathrm{hd}_i M=0$, then the following sequence 
         \begin{equation}\label{eq 3.1}
           M_i\stackrel{(\epsilon(h)M_h)}\lrw \bigoplus_{h:i\to j} M_j \stackrel{(M_{\bar{h}})}\lrw M_i
         \end{equation}
         is exact at the last term, i.e. the map $M_{out}:=(M_{\bar{h}}): \bigoplus\limits_{h:i\to j} M_j\to M_i$ is surjective. On the other hand, the complex $\Hom^{\cdot}(P_{S_i},M)$ of vector spaces is precisely the sequence (\ref{eq 3.1}). Hence, for $M$ with trivial $i$-head, $\Hom^{\cdot}(P_{S_i},M)\otimes P_{S_i}$ is quasi-isomorphic to $(M_i\to \Ker M_{out})\otimes P_{S_i}$, which is equal to $(V_0\oplus V_1[-1])\otimes P_{S_i}$, here $V_0$ (resp. $V_1$) is the $0$-th (resp. $1$-th) cohomology of the complex $M_i\to \Ker M_{out}$. Then $T_{S_i}(M)=\mathrm{cone}( (f_0,f_1):(S_i^{v_0}\oplus S_i^{v_1}[-1])\to M)$, with $f_0=(d_1,\cdots d_{v_0}):S_i^{v_0}\to M$ and $f_1=(c_1,\cdots, c_{v_1}):S_i^{v_1}\to M$, where $d_1,\cdots, d_{v_0}$ forms a basis for $\Hom_{\varLambda_Q}(S_i,M)$ and $c_1\cdots, c_{v_1}$ forms a basis for $\operatorname{Ext}^1_{\varLambda_Q}(S_i,M)$. Denote by $M'$ the extension of $M$ by $S_i^{v_1}$ representing the class $f_1$ in $\operatorname{Ext}^1_{\varLambda_Q}(S_i^{v_1},M)$, then $T_{S_i}(M)=\mathrm{cone}(f_0: S_i^{v_0}\to M')$. Note that $f_0$ is injective and $v_{i_1}+\dim M_i-v_{i_0}=\dim \Ker M_{out}$, it follows that $T_{S_i}(M)$ is quasi-isomorphic to $\Sigma_i(M)$.
   
         For $M\in \mathrm{rep}^n\varLambda_Q$ with $\mathrm{hd}_i M=S_i^v$, notice that there exists a submodule $N$ with trivial $i$-head such that $M$ is an extension of $N$ by $S_i^{v}$. So there is a short exact sequence as follows:
            $$0\lrw N\lrw M\lrw S_i^{v}\lrw 0.$$ 
         Applying $T_{S_i}$, we have a triangle 
             $$T_{S_i}(N)\lrw T_{S_i}(M) \to T_{S_i}(S_i^v)\lrw T_{S_i}(N)[1].$$
         Take $H^{\cdot}$ to the above triangle and note $T_{S_i}=S_i[-1]$ and $T_{S_i}(N)=\Sigma_i(N)$, it follows that $H^0(T_{S_i}(M))\cong H^0(T_{S_i}(N))$. Since $\Sigma_i(M)=\Sigma_i(N)$, we conclude that $H^0(T_{S_i}(M))\cong \Sigma_i(M)$ for any $M\in \mathrm{rep}^n\varLambda_Q$.
        \eremark

    \subsection{Realization of Weyl groups}
        Let $W_Q$ be the Weyl group of the associated Kac-Moody Lie algebra $\ssg_Q$.  Namely,  $W_Q$ is the subgroup of  $\mathrm{Aut}(\Z Q_0)$ generated by simple reflections $s_i$, $i\in Q_0$, where $s_i(\lambda)=\lambda-(\lambda,\alpha_i)_A\alpha_i$ for $\lambda\in \Z Q_0$ and $\{\alpha_i\}$ is the set of simple roots.
   
        Recall that we have an isomorphism $\eta: K_0(\sR_{\varLambda_Q})\to \Z Q_0$, $\hat{S_i}\mapsto \alpha_i$. Let $\boldsymbol{T}_i:K_0(\sR_{\varLambda_Q})\to K_0(\sR_{\varLambda_Q})$ be the induced map of the spherical twist  functor $\bar{T}_{2,S_i}$ of simple $\varLambda_Q$-module $S_i$. Then, we have the following 
      \bproposition\label{prop3.10}
         For any $i\in I$, the following identity holds in $\mathrm{Aut}(\Z Q_0)$:
           $$s_i=\eta \boldsymbol{T}_i \eta^{-1}.$$
        That is, the following diagram commutes
        $$\begin{tikzcd}
        K_0(\sR_{\varLambda_Q}) \arrow[d, "\eta"] \arrow[r, "\boldsymbol{T}_i"]
        & K_0(\sR_{\varLambda_Q}) \arrow[d, "\eta" ] \\
        \Z Q_0 \arrow[r, "s_i"]      & \Z Q_0. \end{tikzcd}$$
      \bproof
       It suffices to check $s_i=\eta \boldsymbol{T}_i \eta^{-1}$ for $\alpha_j$, $j\in Q_0$.  $\eta \boldsymbol{T}_i \eta^{-1}(\alpha_i)=\eta \boldsymbol{T}_i (\hat{S_i})=\eta(-\hat{S_i})=-\alpha_i$. For $j\neq i$, $\eta \boldsymbol{T}_i \eta^{-1}(\alpha_j)=\eta \boldsymbol{T}_i (\hat{S_j})=\eta(\hat{L}_{ji})=\alpha_j+a_{ji}\alpha_i$. Combining with $s_i(\alpha_i)=-\alpha_i$ and $s_i(\alpha_j)=\alpha_j-(\alpha_j,\alpha_i)_A\alpha_i=\alpha_j+a_{ji}\alpha_i$, where $-a_{ji}=<\hat{S_j},\hat{S_i}>_{\varLambda_Q}=(\hat{S_j},\hat{S_i})_A$, we complete the proof.
      \eproof
      \eproposition
      In the sequel, we will regard $T_i$ as $s_i$ in $W_Q$.

      Let $\mathrm{Br}(\sR_{\varLambda_Q})$ be the subgroup of the autoequivalence group $\mathrm{Auteq}(\sR_{\varLambda_Q})$ generated by $\bar{T}_{2,S_i}$, $i\in Q_0$. It can be seen that the map $\epsilon:\mathrm{Br}(\sR_{\varLambda_Q})\to W_Q$ taking any autoequivalence to its restriction on $K_0(\sR_{\varLambda_Q})$ is a surjective group homomorphism following from Proposition \ref{prop3.10}.  On the other hand, there is a surjective group homomorphism $\xi:\mathrm{Br}(\sR_{\varLambda_Q})\to \mathrm{Aut}_0(\ssg(\sR_{\varLambda_Q}))$, $\bar{T}_{2,S_i}\mapsto \Phi_i$, where $\mathrm{Aut}_0(\ssg(\sR_{\varLambda_Q}))$ is the subgroup of $\mathrm{Aut}(\ssg(\sR_{\varLambda_Q}))$ generated by $\Phi_i$.
   
      \btheorem\label{thm3.11}
      There exists a group homomorphism
         $$\epsilon':\mathrm{Aut}_0(\ssg(\sR_{\varLambda_Q}))\to W_Q,\ \ \Phi_i\mapsto s_i.$$
      Moreover, $\epsilon'$ is surjective.
      \bproof
       Let $\Ker\xi$ be the kernel of $\xi:\mathrm{Br}(\sR_{\varLambda_Q})\to \mathrm{Aut}_0(\ssg(\sR_{\varLambda_Q}))$, we claim that for any $w\in \Ker\xi$, $\epsilon(w)=0$. Note that for each element $w\in \mathrm{Br}(\sR_{\varLambda_Q})$, $w$ is a product of $\bar{T}_{2,S_i}$, write $w=\bar{T}_{2,S_{i_1}}\bar{T}_{2,S_{i_2}}\cdots\bar{T}_{2,S_{i_r}}$. Then for $w\in \Ker\xi$, by the definition of $\Phi_i$, we have $w(X)\cong X$ for $X\in \mathrm{ind}\sR_{\varLambda_Q}$, it follows that $\boldsymbol{T}_{i_1}\boldsymbol{T}_{i_2}\cdots \boldsymbol{T}_{i_r}([S_i])=[S_i]$ for any $i\in Q_0$. Hence $\epsilon(w)(\alpha_i)=\eta \boldsymbol{T}_{i_1}\boldsymbol{T}_{i_2}\cdots \boldsymbol{T}_{i_r}\eta^{-1}(\alpha_i)=\alpha_i$ for $i\in Q_0$, which means that $\epsilon(w)=1$. Then, there exists a group homomorphism $\epsilon':\mathrm{Aut}_0(\ssg(\sR_{\varLambda_Q}))\to W_Q$ such that the following diagram commutes.
       $$\begin{tikzcd}
        \mathrm{Br}(\sR_{\varLambda_Q}) \arrow[d, "\epsilon"] \arrow[r, "\xi"]
        & \mathrm{Aut}_0(\ssg(\sR_{\varLambda_Q})) \arrow[dl, "\epsilon'" ] \\
        W_Q
      \end{tikzcd}$$
      \eproof
      \etheorem
   
      \bremark
       As shown in \cite[Section 5]{Baumann2012a}, the reflection functor $\Sigma_i$ is the simple reflection $s_i\in W_Q$ at the crystal level in some sense. However the reflection functor $\Sigma_i$ may not be well defined at the weight $\alpha_i$. Because $s_i(\alpha_i)=-\alpha_i$ which is not a positive root anymore. In order to overcome it, we extend the reflection functor $\Sigma_i$ on $\mathrm{rep}^n\varLambda_Q$ to the spherical twist functor $\bar{T}_i$ on the triangulated hull $\sR_{\varLambda_Q}$ of $D^b(\mathrm{rep}^n\varLambda_Q)$, and reach the result Theorem \ref{thm3.11}. 
      \eremark
   
\section{Relations between \texorpdfstring{$\mathfrak{g}_0
      (\mathcal{R}_A)_{(q-1)}$}{Lg} and \texorpdfstring{$\mathfrak{g}_0(\mathcal{R}_{\varLambda_Q})_{(q-1)}$}{Lg}}
 \subsection{Relations between \texorpdfstring{$\mathfrak{n}^+_0
      (\mathcal{R}_A)_{(q-1)}$}{Lg} and \texorpdfstring{$\mathfrak{n}^+_0(\mathcal{R}_{\varLambda_Q})_{(q-1)}$}{Lg}}
          Let $k$ be the finite field with $|k|=q$. Recall that $\ssg_0(\sR_A)$ (resp. $\ssg_0(\sR_{\varLambda_Q})$) is a subalgebra of $\ssg(\sR_A)_{(q-1)}$ (resp. $\ssg(\sR_{\varLambda_Q})_{(q-1)}$) generated by $u_{S_j}$, $u_{S_j[1]}$, $h_j$,  $j\in Q_0$ (resp. $\hat{u}_{S_j}$, $\hat{u}_{S_j[1]}$, $h_j$, $j\in Q_0$), and $\ssn^+_0(\sR_A)$ (resp. $\ssn^+_0(\sR_{\varLambda_Q})$) generated by $u_{S_j}$, $j\in Q_0$ (resp. $\hat{u}_{S_j}$, $j\in Q_0$).
      
          For $M,N,L\in \mathrm{rep}^n\varLambda_Q$, denote by $\Hom_{\varLambda_Q}(N,L)_M$ the subset of $\Hom_{\varLambda_Q}(N,L)$ consisting of all injective homomorphisms whose $\textbf{cokernel}$ is isomorphic to $M$. Then for any $f_0\in \Hom_{\varLambda_Q}(N,L)_M$, we have a triangle in $D^b(\mathrm{rep}^n\varLambda_Q)$ as follows:
          $$N\stackrel{f_0}\lrw L\lrw M \lrw N[1].$$
         Hence its image under the functor $\pi: D^b(\mathrm{rep}^n\varLambda_Q)\to \sR_{\varLambda_Q}$ is also a triangle in $\sR_{\varLambda_Q}$.
          So $\Hom_{\varLambda_Q}(N,L)_M\subset \Hom_{\sR_{\varLambda_Q}}(N,L)_M$ by regarding $f_0$ as $(f_0,0)$, because $\Hom_{\sR_{\varLambda_Q}}(N,L)\cong \Hom_{\varLambda_Q}(N,L)\oplus \operatorname{Ext}^2_{\varLambda_Q}(N,L)$.
      
         Set
              $$\Hom_{\sR_{\varLambda_Q}}(N,L)_M:=\{f\in \Hom_{\sR_{\varLambda_Q}}(N,L)\ |\ \mathrm{cone}f\cong M\}.$$
      
       \bexample
         Take the acyclic quiver $Q$ to be the Kronecker quiver
         $$\begin{tikzcd} 1\arrow[r,"\beta_1",shift left=0.5ex] \arrow[r,"\beta_2"',shift right=0.5ex] &2\end{tikzcd}$$
          Let us compute $[S_1][S_2]$ in $\mathrm{rep}^n\varLambda_Q$ and in $\sR_{\varLambda_Q}$ respectively.
          $$[S_1][S_2]=\sum\limits_{\lambda\in k}[E_{\lambda}]+[E']+[S_1\oplus S_2]  $$
          in the Hall algebra of $\mathrm{rep}^n\varLambda_Q$, where $E_{\lambda}$ and $E'$ are as follows:
          $$\begin{tikzcd}E_{\lambda}:\  k\arrow[r,"1",shift left=0.5ex] \arrow[r,"\lambda"',shift right=0.5ex] &k\end{tikzcd}\qquad\ and\qquad
          \begin{tikzcd}E':\  k\arrow[r,"0",shift left=0.5ex] \arrow[r,"1"',shift right=0.5ex] &k\end{tikzcd}.$$
           On the other hand, $S_2\lrw L\lrw S_1\stackrel{w}\lrw S_2[1]$ is a triangle in $\sR_{\varLambda_Q}$ if and only if $L[1]\cong \mathrm{cone}w$. Note $w\in \Hom_{\sR_{\varLambda_Q}}(S_1,S_2[1])=\operatorname{Ext}^1_{\varLambda_Q}(S_1,S_2)$, then $L\in D^b(\mathrm{rep}^n \varLambda_Q)$, actually $L$ is an extension of $S_1$ by $S_2$. $L\cong S_1\oplus S_2$ iff $w=0$, then
             $$F_{S_1,S_2}^{S_1\oplus S_2}=\frac{\Hom_{\sR_{\varLambda_Q}}(S_2,S_1\oplus S_2)_{S_1}}{\mathrm{Aut}_{\sR_{\varLambda_Q}}(S_2)}=\frac{\Hom_{\sR_{\varLambda_Q}}(S_2,S_2)_0}{\mathrm{Aut}_{\sR_{\varLambda_Q}}(S_2)}=1.$$
        If $L\cong E_{\lambda}$, Applying $\Hom_{\varLambda_Q}(S_2,-)$ to the short exact sequence $S_2\hookrightarrow E_{\lambda}\twoheadrightarrow S_1$, we have a long exact sequence
        $$0\to \operatorname{Ext}^1_{\varLambda_Q}(S_2,E_{\lambda})\to  \operatorname{Ext}^1_{\varLambda_Q}(S_2,S_1)\overset{w_{\lambda*}}\to  \operatorname{Ext}^2_{\varLambda_Q}(S_2,S_2)\to  \operatorname{Ext}^2_{\varLambda_Q}(S_2,E_{\lambda})\to 0$$
        Since $w_{\lambda}\neq 0$, we have that $w_{\lambda*}$ is an isomorphism. Hence
        $$\Hom_{\sR_{\varLambda_Q}}(S_2,E_{\lambda})=\Hom_{\varLambda_Q}(S_2,E_{\lambda}).$$
         Then
           $$F_{S_1,S_2}^{E_{\lambda}}=\frac{\Hom_{\sR_{\varLambda_Q}}(S_2,E_{\lambda})_{S_1}}{\mathrm{Aut}_{\sR_{\varLambda_Q}}(S_2)}=\frac{\Hom_{{\varLambda_Q}}(S_2,E_{\lambda})_0}{\mathrm{Aut}_{\varLambda_Q}(S_2)}=1.$$
        Similarly, we have that $F_{S_1,S_2}^{E'}=1$. So
             $$[S_1][S_2]=\sum\limits_{\lambda\in k}[E_{\lambda}]+[E']+[S_1\oplus S_2]  $$
           in the Hall algebra of $\sR_{\varLambda_Q}$.
       \eexample
       Let $M$, $N$, be indecomposable $\varLambda_Q$-modules, and $L$ an extension of $M$ by $N$. We may not be able to compute $\Hom_{\sR_{\varLambda_Q}}(N,L)$. However we still can show that $F_{M,N}^{L}(\mathrm{rep}^n\varLambda_Q)\equiv F_{M,N}^{L}(\sR_{\varLambda_Q})\pmod{(q-1)}$, where $L$ is an indecomposable $\varLambda_Q$-module. Here $F_{M,N}^L(\sA)$ means the cardinality of $\frac{\Hom_{\sA}(N,L)_M}{\mathrm{Aut}_{\sA}(N)}$ in the given category $\sA$.
      
       \blemma\label{lem2}
          Let $M$, $N\in \mathrm{rep}^n\varLambda_Q$ be indecomposable, and $L\in \sR_{\varLambda_Q}$ indecomposable such that
          $$N\stackrel{f}\lrw L\lrw M\lrw N[1]$$
          is a triangle in $\sR_{\varLambda_Q}$, then $L\in \mathrm{rep}^n\varLambda_Q$, and
            $$\Hom_{\sR_{\varLambda_Q}}(N,L)_M=\{(f_0,f_1)|\mathrm{cone}f_0\cong M\}=\Hom_{\varLambda_Q}(N,L)_M\oplus \operatorname{Ext}^2_{\varLambda}(N,L).$$
      
        \bproof
          Note that if $N\stackrel{f}\to L\stackrel{g}\to M\stackrel{h}\to N[1]$ is a triangle in $\sR_{\varLambda_Q}$ with $M,N\in \mathrm{rep}^n\varLambda_Q$, then $h\in \Hom_{\sR_{\varLambda_Q}}(M,N[1])=\operatorname{Ext}^1_{\varLambda_Q}(M,N)$, since $gl.dim \varLambda_Q=2$, it follows that $L[1]\cong \mathrm{cone}h\in D^b(\mathrm{rep}^n\varLambda_Q)$ must be an extension of $M$ by $N$. Moreover, for any element $(f_0,f_1)$ of $\Hom_{\sR_{\varLambda_Q}}(N,L)_M$, $f_0\neq 0$. Indeed, if $f_0=0$ then there is a triangle in $\sR_{\varLambda_Q}$,
             $$N\stackrel{f_1}\lrw L\stackrel{g}\lrw M\stackrel{h}\lrw N[1].$$
          But the functor $\pi: D^b(\mathrm{rep}^n\varLambda_Q)\to \sR_{\varLambda_Q}$ is a triangle functor and $f_1\in \Hom_{D^b(\varLambda_Q)}(N,L[2])$, the above triangle must be the image of the triangle $N\stackrel{f_1}\lrw M \lrw \mathrm{cone}f_1\lrw N[1]$ in $D^b(\mathrm{rep}^n\varLambda_Q)$. Hence $M[2l]\cong c\mathrm{cone}f_1$ in $D^b(\mathrm{rep}^n\varLambda_Q)$ for some integer $l$, which implies that $g=0$ or $h=0$. So $N\cong M[-1]\oplus L[2]$ or $L[2]\cong N\oplus M[2l]$, this is a contradiction.
      
          Let $f=(f_0,f_1)\in \Hom_{\sR_{\varLambda_Q}}(N,L)_M$, then we have a triangle in $\sR_{\varLambda_Q}$
             $$N\stackrel{(f_0,f_1)}\lrw L\stackrel{(g_0,g_1)}\lrw M\stackrel{h}\lrw N[1].$$
          then $(f_0,f_1)h[-1]=(f_0h[1],f_1h[-1])=0$, this implies that $f_0h[-1]=0$. Hence there exists $v=(v_0,v_1)\in \Hom_{\sR_{\varLambda_Q}}(L,L)$ such that $vf=f_0$. So we have a morphism of triangles
          $$\begin{tikzcd}
            N\arrow[d, "id"] \arrow[r, "{(f_0,f_1)}"]
            &L \arrow[r, "g" ]\arrow[d,"{(v_0,v_1)}"]   &M\arrow[d,dotted]\arrow[r]   &N[1]\arrow[d]\\
            N\arrow[r,"f_0"] &L\arrow[r]  &cone(f_0)\arrow[r]  &N[1].
            \end{tikzcd}$$
          If $v_0\in \operatorname{End}_{\varLambda_Q}(L)$ is a unit, then the two triangles above are isomorphic and we get that $\mathrm{cone}f_0\cong M$.
      
          If $v_0\in \operatorname{End}_{\varLambda_Q}(L)$ is not a unit, then $v_0^r=0$ for some positive integer $r$ because $\operatorname{End}_{\varLambda_Q}(L)$ is a local ring. By direct computation, $(v_0,v_1)^{2r}=0$. Therefore, we have a morphism as follows between the two triangles above
          $$\begin{tikzcd}
            N\arrow[d, "id"] \arrow[r, "{(f_0,f_1)}"]
            &L \arrow[r, "g" ]\arrow[d,"{(v_0,v_1)^{2r}}"]   &M\arrow[d,dotted]\arrow[r]   &N[1]\arrow[d]\\
            N\arrow[r,"f_0"] &L\arrow[r]  &cone(f_0)\arrow[r]  &N[1].
           \end{tikzcd}$$
          Then $f_0=v^{2r}f=0$, this contradicts to $f_0\neq 0$.
      
          Therefore, we have shown that $\Hom_{\sR_{\varLambda_Q}}(N,L)_X\subset \{(f_0,f_1)|\mathrm{cone}f_0\cong M\}$.
      
          On the other hand, let $f=(f_0,f_1)\in \Hom_{\sR_{\varLambda_Q}}(N,L)$ such that $\mathrm{cone}f_0\cong M$. Consider the triangle
          $$N\stackrel{(f_0,f_1)}\lrw L\lrw \mathrm{cone}f \lrw N[1].$$
          Applying the functor $\Hom_{\varLambda_Q}(-,L)$ to the short exact sequence $0\to N\stackrel{f_0}\to L\stackrel{g_0}\to M\to 0$, we get a long exact sequence
            $$\cdots \to \operatorname{Ext}_{\varLambda_Q}^2(L,L) \to \operatorname{Ext}_{\varLambda_Q}^2(Y,L)\to 0,$$
           since $\operatorname{Ext}^3_{\varLambda_Q}(X,L)=0$. Therefore, for any $f_1\in \operatorname{Ext}^2_{\varLambda_Q}(Y,L)$, there exists $v_1\in \operatorname{Ext}^2(L,L)$ such that $f_1=v_1f_0$. Note that $L$ is indecomposable and $\operatorname{Ext}^2_{\varLambda_Q}(L,L)\subset \mathrm{rad} (\Hom_{\sR_{\varLambda_Q}}(L,L))$, so $(id,v_1)$ is an isomorphism in $\Hom_{\sR_{\varLambda_Q}}(L,L)$. Then the following two triangles are isomorphic
           $$\begin{tikzcd}
            N\arrow[d, "id"] \arrow[r, "f_0"]
            &L \arrow[r, "g" ]\arrow[d,"{(id,v_1)}"]   &M\arrow[d,dotted]\arrow[r]   &N[1]\arrow[d]\\
            N\arrow[r,"{(f_0,f_1)}"] &L\arrow[r]  &cone(f_0)\arrow[r]  &N[1].
          \end{tikzcd}$$
          which implies $\mathrm{cone}f\cong M$ and $f\in \Hom_{\sR_{\varLambda_Q}}(N,L)_M$. The proof is completed.
        \eproof
        \elemma
      
          The automorphism group $\mathrm{Aut}_{\varLambda_Q}(N)$ can be naturally embedded into the automorphism group $\mathrm{Aut}_{\sR_{\varLambda_Q}}(N)$ by $f_0\mapsto (f_0,0)$, and the extension group $\operatorname{Ext}^2_{\varLambda_Q}(N,N)$ also can be regarded as a subgroup of $\mathrm{Aut}_{\sR_{\varLambda_Q}}(N)$ by $f_1\mapsto (1,f_1)$. Furthermore, $\mathrm{Aut}_{\varLambda_Q}(N)$ acts on $\mathrm{Ext}^2_{\varLambda_Q}(N,N)$ by composition, and for any $(f_0,f_1)$ in $\mathrm{Aut}_{\sR_{\varLambda_Q}}(N)$, $(f_0,f_1)=(f_0,0)(1,f_0^{-1}f_1)$ $\in \mathrm{Aut}_{\varLambda_Q}(N)\ltimes \operatorname{Ext}^2_{\varLambda_Q}(N,N)$. Hence we have that
             $$\mathrm{Aut}_{\sR_{\varLambda_Q}}(N)= \mathrm{Aut}_{\varLambda_Q}(N)\ltimes \operatorname{Ext}^2_{\varLambda_Q}(N,N).$$
      
        \blemma\label{lem3}
          For $M,N,L$ indecomposable in $\mathrm{rep}^n\varLambda_Q$, we have that
           $$|\operatorname{Ext}^2_{\varLambda_Q}(N,N)|\cdot|\frac{\Hom_{\sR_{\varLambda_Q}}(N,L)_M}{\mathrm{Aut}_{\sR_{\varLambda_Q}}(N)}|=|\operatorname{Ext}^2_{\varLambda_Q}(N,L)|\cdot|\Hom_{\varLambda_Q}(N,L)^*_M| .$$
        \elemma
        \bproof
          By Lemma \ref{lem2} we know that $\Hom_{\sR_{\varLambda_Q}}(N,L)_M\cong \Hom_{\varLambda_Q}(N,L)_M$ $\oplus  \operatorname{Ext}^2_{\varLambda_Q}(N,L)$, it follows that $f_0$ is injective for any $f=(f_0,f_1)\in \Hom_{\sR_{\varLambda_Q}}(N,L)_M$. Thus, the action of $\mathrm{Aut}_{\sR_{\varLambda_Q}}(N)$ on $\Hom_{\sR_{\varLambda_Q}}(N,L)_M$ is free. Indeed, if $(f_0,f_1)(u,v)=(f_0,f_1)$, then $f_0u=f_0$ and $f_0v+f_1u=f_1$, which implies $u=id_N$ and $v=0$. Naturally, as a subgroup of $\mathrm{Aut}_{\sR_{\varLambda_Q}}(N)$, $\mathrm{Aut}_{\varLambda_Q}{N}$ acts on $\Hom_{\sR_{\varLambda_Q}}(N,L)_{M}$ freely. Then we have a surjection
                $$E:= \frac{\Hom_{\varLambda_Q}(N,L)_M\oplus \operatorname{Ext}^2_{\varLambda_Q}(N,L)}{\mathrm{Aut}_{\varLambda_Q}(N)}\longrightarrow \frac{\Hom_{\varLambda_Q}(N,L)_M}{\mathrm{Aut}_{\varLambda_Q}(N)}$$
           with fiber $\operatorname{Ext}^2_{\varLambda_Q}(N,L)$.
      
          On the other hand, since $\mathrm{Aut}_{\sR_{\varLambda_Q}}(N)= \mathrm{Aut}_{\varLambda_Q}(N)\ltimes \operatorname{Ext}^2_{\varLambda_Q}(N,N)$, we have another surjection
                $$E= \frac{\Hom_{\sR_{\varLambda_Q}}(N,L)_M}{\mathrm{Aut}_{\varLambda_Q}(N)}\longrightarrow \frac{\Hom_{\sR_{\varLambda_Q}}(N,L)_M}{\mathrm{Aut}_{\varLambda_Q}(N)\ltimes \operatorname{Ext}^2_{\varLambda_Q}(N,N)}$$
          with fiber $\operatorname{Ext}^2_{\varLambda_Q}(N,N)$. We complete the proof.
         \eproof
      
         Since the integral Hall algebra $\overline{\sH(\varLambda_Q)}$ of the preprojective algebra is an associative algebra, we can define a Lie algebra $\ssn^{+}(\varLambda_Q)$ to be the free abelian group with the basis $\{[M]|M\in rep^n\varLambda_Q\}$, and the Lie bracket is given by $[[M],[N]]=[M][N]-[N][M]$. Also, we can consider the Lie subalgebra $\ssn^+_0(\varLambda_Q)$ of $\ssn^+(\varLambda_Q)$ generated by $[S_i]$ for $i\in I$.
      
         \bproposition\label{prop4.3}
          There is an isomorphism of Lie algebras:
                 $$\phi_0:\ssn^+_0(\sR_{\varLambda_Q})\lrw \ssn^+_0(\varLambda_Q),\ \hat{u}_{S_i}\mapsto [S_i].$$
          \bproof
           By Lemma \ref{lem3} and Lemma \ref{lem2}, we have an algebra homomorphism from $\phi:\ssn^+(\sR_{\varLambda_Q})\lrw \ssn^+(\varLambda_Q)$, $\hat{u}_M\mapsto [M]$, which is obviously a bijection. Thus, $\phi$ maps $\ssn^+_0(\sR_{\varLambda_Q})$ to $\ssn^+_0(\varLambda_Q)$ bijectively, which is exactly $\phi_0$.
          \eproof
         \eproposition
      
         Let us recall Lusztig's categorification of the positive half part of the enveloping algebra $U^+$ of the associated Kac-Moody Lie algebra $\ssg_Q$ using constructible functions \cite{Lusztig2000}. Fix a dimension vector $\boldsymbol{v}\in \N Q_0$, $G_{\boldsymbol{v}}:=\prod_{i\in I}GL(v_i;\C)$. Consider the representation space $\boldsymbol{E}_{\boldsymbol{v}}=\prod_{h\in Q_1}\Hom(\C^{v_{s(h)}},\C^{v_{t(h)}})$ consisting of all $\C Q$-modules with dimension vector $\boldsymbol{v}$. Note $\boldsymbol{E}_{\boldsymbol{v}}$ is an affine variety over complex field $\C$. Denote by $M(\boldsymbol{E}_{\boldsymbol{v}})$ the $\Q$-vector space of all constructible $G_{\boldsymbol{v}}$-equivariant functions $f:\varLambda_{\boldsymbol{v}} \to \Q$ (A constructible function is a finite sum of constant functions on a constructible subset of $\boldsymbol{E}_{\boldsymbol{v}}$ and $f$ is said to be $G_{\boldsymbol{v}}$-equivariant if $f$ is constant on $G_{\boldsymbol{v}}$-orbits ). Set $M(\boldsymbol{E}):=\bigoplus_{\boldsymbol{v}\in \N I} M(\boldsymbol{E}_{\boldsymbol{v}})$, it has a standard convolution product $*$. Namely, for $f\in M(\boldsymbol{E}_{\boldsymbol{v}'})$, $g\in M(\boldsymbol{E}_{\boldsymbol{v}''})$, set $\boldsymbol{v}=\boldsymbol{v}'+\boldsymbol{v}''$, define
              $$f*g(x):=\sum_{x'\in \boldsymbol{E}_{\boldsymbol{v}'},x''\in \boldsymbol{E}_{\boldsymbol{v}''}}\chi(V(x',x'';x))f(x')g(x'')$$
         where $V(x',x'';x)=\{L\subset x|L\cong x'', x/L\cong x'\}$, which is a locally closed subset of $\boldsymbol{E}_{\boldsymbol{v}}$, and $\chi(V)$ is the Euler character of a variety $V$ (see \cite [Lemma 3.3]{Kasjan2010}).
      
          Consider the $\Q$-linear subalgebra $M_0(\boldsymbol{E})$ of $M(\boldsymbol{E})$ generated by $f_i=1_{\sO(S_i)}$($i\in Q_0$), where $\sO(S_i)$ is the $G_v$-orbit of $S_i$. By \cite[Proposition 10.20]{Lusztig2000} , there is an isomorphism
              $$\kappa : U^+\lrw M_0(\boldsymbol{E}),\ \ e_i\mapsto f_i,$$
          where $U^+$ is the  $\Q$-algebra $U^+$ generated by $e_i$, $i\in Q_0$ and subject to the Serre relations
                  \begin{equation}\label{eqn4.1}
                  \sum_{l=0}^{N}(-1)^l \left(\begin{smallmatrix}N \\ l\end{smallmatrix}\right) e_i^{N-l}e_j e_i^l=0,\end{equation}
         for $i\neq j$, where $N=1+a_{ij}+a_{ji}$, $a_{ij}=\dim \operatorname{Ext}^1_A(S_i,S_j)$.
      
        Let $U^+_{\Z}$ be the $\Z$-linear subalgebra of $M_0(\boldsymbol{E})$ generated by $f_i$, $i\in Q_0$. Note Serre relations (\ref{eqn4.1}) are $\Z$-coefficients, hence $U^+_{\Z}$ is a $\Z$-algebra generated by $f_i$, $i\in Q_0$ and subject to
               \begin{equation*}
               \sum_{l=0}^{N}(-1)^l \left(\begin{smallmatrix}N \\ l\end{smallmatrix}\right) f_i^{N-l}f_j f_i^l=0,\end{equation*}
          for $i\neq j$, where $N=1+a_{ij}+a_{ji}$.
          Finally, we define a $\Z/(q-1)$-algebra by setting
             $$\overline{U^+_{\Z}}:=U^+_{\Z}\otimes_{\Z} \Z/(q-1).$$
          It is clear that $\overline{U^+_{\Z}}$ is a $\Z/(q-1)$-algebra with generators $\bar{f_i}$, $i\in Q_0$, subject to
          \begin{equation*}
            \sum_{l=0}^{N}(-1)^l \overline{\left(\begin{smallmatrix}N \\ l\end{smallmatrix}\right)} \bar{f_i}^{N-l}\bar{f_j} \bar{f_i}^l=0,\end{equation*}
            for $i\neq j$, where $N=1+a_{ij}+a_{ji}$.
      
         We introduce some notations. 
            $$\left(\begin{smallmatrix}m \\ l\end{smallmatrix}\right):=\frac{m!}{l!(m-l)!},\qquad [l]_q:=\frac{q^l-1}{q-1},\qquad [\begin{smallmatrix}m \\ l\end{smallmatrix}]_q:=\frac{[m]_q!}{[l]_q![m-l]_q!}.$$
         Fix a dimension vector $\boldsymbol{v}\in\N Q_0$. Set
           $$\aleph_{\boldsymbol{v}}:=\{\lambda=l_1i_1+l_2i_2+\cdots l_ti_t|i_j\in I,\sum_{j=1}^t l_ji_j=\boldsymbol{v}, l_j\in \N\}$$
        and $\aleph:=\bigcup\limits_{\boldsymbol{v}\in \N I} \aleph_{\boldsymbol{v}}$. For  $\lambda=l_1i_1+l_2i_2+\cdots l_ti_t\in \aleph_{\boldsymbol{v}}$ and $M\in \mathrm{mod}A$, to compute $[S_{i_1}]^{l_1}$ $[S_{i_2}]^{l_2}\cdots [S_{i_t}]^{l_t}\in C\sH(A)$, we first compute $S_{\lambda}:=[S_{i_1}^{\oplus l_1}][S_{i_2}^{\oplus l_2}]\cdots [S_{i_t}^{\oplus l_t}]$. Denote by $F_{\lambda}^M$ the number of times of $[M]$ in  $[S_{i_1}^{\oplus l_1}][S_{i_2}^{\oplus l_2}]\cdots [S_{i_t}^{\oplus l_t}]$, then
          $$[S_{i_1}]^{l_1}[S_{i_2}]^{l_2}\cdots [S_{i_t}]^{l_t}= \prod_{j=1}^t [l_j]_q!\sum_{[M]} F_{\lambda}^M [M].$$
        Let
          $$V(\lambda;M)_q:=\{L_t\subset \cdots \subset L_1=M|L_j/L_{j+1}\cong S_{i_j}^{\oplus l_j} \}.$$
        Then  the cardinality $|V(\lambda;M)_q|$ is exactly the number $F_{\lambda}^M$ from the definition of Hall product.
      
        On the other hand, for $f_{i_1}^{l_1}f_{i_2}^{l_2}\cdots f_{i_t}^{l_t}\in M_0(\boldsymbol{E}_v)$, note $f_{i}^l=\prod_{j=1}^{l-1}\chi(\mathrm{Gr}(j,j+1))$ $1_{\sO(S_i^{\oplus l})}$= $(l!) 1_{\sO(S_i^{\oplus l})}$, here $\mathrm{Gr}(m,n)$ means the Grassmannian variety of $m$-dimensional vector subspaces of the vector space $\C^n$.  Set
                $$V(\lambda;x)_{\C}:=\{x_t\subset \cdots \subset x_1=x|x_j\in \mathrm{mod}\C Q,\ x_j/x_{j+1}\cong S_{i_j}^{\oplus l_j} \}.$$
        Then $f_{i_1}^{l_1}f_{i_2}^{l_2}\cdots f_{i_t}^{l_t}=\prod_{j=1}^t (l_j)!\sum\limits_{[M]} \chi(V(\lambda;M)_{\C})1_{\sO(M)}$.
      
        \blemma\label{lem5}
         Regard $q=|k|$ as an indeterminate, then for $\lambda\in \aleph_v$ and $M\in modA$, we have that $F_{\lambda}^{M}\in \Z[q]$.
         \bproof
           We compute $F_{\lambda}^M=|V(\lambda;M)|$ by induction on $|v|$ for $v\in \N I$, $\lambda\in \aleph_v$, and $M\in modA$. If $|v|=1$, then $S_{\lambda}=S_{i}$ for some $i\in I$. It is obvious that $F_{\lambda}^{M}$ is a constant function. Assume the statements holds for $|v|<m$, for $|v|=m$ and $\lambda=l_1i_1+\cdots l_ti_t\in \aleph_v$, consider all submodules $N$ of $M$ such that $N\cong S_{i_t}^{\oplus l_t}$, define an equivalent relation by $N\sim N'$ iff $M/N\cong M/N'$. Denote the class of $N$ by $\hat{N}$ and the set of all equivalent classes by $\Omega$. Then
                    $$F_{\lambda}^M=\sum_{\hat{N}\in \Omega} [\begin{smallmatrix}d_M-a_N\\l_t-a_N \end{smallmatrix}]_qF_{\lambda-l_ti_t}^{M/N},$$
          where $d_M=\dim \bigcap_{s(h)=i_t} \Ker h$, $a_N=\dim N\cap \sum_{t(h)=i_t} \mathrm{im} h$. Note $\lambda-l_ti_t\in \aleph_{v'}$ with $|v'|<|v|$, by induction, $F_{\lambda-l_ti_t}^{M/N}\in \Z[q]$. It is well-known that $[\begin{smallmatrix}d_M-a_N\\l_t-a_N \end{smallmatrix}]_q\in \Z[q]$. Hence, $F_{\lambda}^M\in \Z[q]$.
         \eproof
        \elemma
      
        Due to Lemma \ref{lem5}, the variety $V(\lambda;M)_{q}$ can be counted by a polynomial, we have $\chi_c(V(\lambda;M)_{\C})=F_{\lambda}^M(1)$ for any $\lambda\in \aleph_v$ and $M\in \mathrm{mod}A$, which can be deduced from \cite[Proposition 6.1]{Reineke2006} or \cite[Lemma 8.1]{Reineke2008}. 
        \blemma[ \cite{Reineke2008}]
         Suppose there exists a rational function $P_X(t) \in \Q(t)$ such that $P_X(|k|) = |X(k)|$ for almost all reductions $k$ of the ring $R$ over which $X$ is defined. Then $P_X(t) \in \Z[t]$ is actually a polynomial with integer coefficients, and $P_X(1) =\chi_c(X)$, the Euler characteristic of $X$ in singular cohomology with compact support.
        \elemma
        Note $V(\lambda;M)_{\C}$ is a locally closed subset defined over an algebraically closed field $\C$, the compactly supported Euler characteristic $\chi_c(V(\lambda;M)_{\C})$ is equal to the Euler characteristic $\chi(V(\lambda;M)_{\C})$. For details, see page 141 in \cite{Fulton1993}.
      
        Set $\overline{C\sH(A)}:=C\sH(A)/ (C\sH(A)\cap (q-1)\sH(A))$, $\overline{C\sH(\varLambda_Q)}:=C\sH(\varLambda_Q)/ (C\sH(\varLambda_Q)\cap (q-1)\sH(\varLambda_Q))$.
      
        \blemma\label{lem4}
          The morphism $\kappa_1: \overline{U^+_{\Z}}\lrw \overline{C\sH(A)}$ sending $\bar{f_i}$ to $[S_i]$ is an isomorphism of $\Z/(q-1)$-algebras.
        \bproof
          First, we should show that $\kappa_1$ is well-defined. Namely, $\overline{C\sH(A)}$ admits the Serre relations
             $$\sum_{l=0}^{N}(-1)^l \left(\begin{smallmatrix}N \\ l\end{smallmatrix}\right) [S_i]^{N-l}[S_j] [S_i]^{l}=0 \pmod{(q-1)\sH(A)},$$
          for $i\neq j$ and $a_{ij}+a_{ji}\neq 0$, where $N=1+a_{ij}+a_{ji}$.
      
          Fix $i\neq j\in I$,  denote by $a_{M,l}$ the multiple of $[M]$ in $[S_i]^{N-l}[S_j] [S_i]^{l}$ for $M\in \mathrm{mod}A$. Note that $Q$ is an acyclic quiver, then either $a_{ij}=0$ and $a_{ji}\neq 0$, or $a_{ij}\neq 0$ and $a_{ji}=0$. In the first case, say there are $a_{ji}$ arrows $\bar{h}$ from $j\to i$. $a_{M,l}\neq 0$ iff $V_M:=\sum_{\bar{h}} \mathrm{im} \bar{h}\subset S_i^{\oplus l}$. If so, $a_{M,l}=[l]_q! [N-l]_q! [\begin{smallmatrix}  N-\bar{d_M}\\l-\bar{d_M} \end{smallmatrix}]_q$, here $\bar{d}_M=\dim V_M\leq a_{ji}<N=1+a_{ji}+a_{ij}$. Then, we have that
                $$\begin{aligned}
                &\sum_{l=0}^{N}(-1)^l \left(\begin{smallmatrix}N \\ l\end{smallmatrix}\right) [S_i]^{N-l}[S_j] [S_i]^{l} \\
                & = \sum_{[M]}\sum_{l=0}^{N}(-1)^l \left(\begin{smallmatrix}N \\ l\end{smallmatrix}\right)[l]_q! [N-l]_q! [\begin{smallmatrix}  N-\bar{d}_M\\l-\bar{d}_M \end{smallmatrix}]_q [M]\\
                &\equiv \sum_{[M]} \sum_{l=0}^{N}(-1)^l N! \left(\begin{smallmatrix}  N-\bar{d}_M\\l-\bar{d}_M \end{smallmatrix}\right)[M] \pmod{(q-1)\sH(A)}\\
                &\equiv 0\pmod{(q-1)\sH(A)}.
                \end{aligned}$$
      
          In the second case, say there are $a_{ij}$ arrows ${h_r}$ from $i$ to $j$. Then $a_{M,l}\neq 0$ iff $S_i^{\oplus I}\subset \bigcap_{h_r} \Ker h_r=:W_M$. If so, $a_{M,l}=[l]_q! [N-l]_q! [\begin{smallmatrix}  d_M\\l \end{smallmatrix}]_q$, here $d_M=\dim W_M\geq 1$ for $\dim \Ker h_r\geq a_{ij}$, $1\leq r\leq a_{ij}$. Thus, we have that
              $$\begin{aligned}
               &\sum_{l=0}^{N}(-1)^l \left(\begin{smallmatrix}N \\ l\end{smallmatrix}\right) S_i^{N-l}S_j S_i^{l} \\
               &= \sum_{[M]} \sum_{l=0}^{N}(-1)^l \left(\begin{smallmatrix}N \\ l\end{smallmatrix}\right)[l]_q! [N-l]_q! [\begin{smallmatrix}  d_M\\l \end{smallmatrix}]_q [M]\\
               &\equiv \sum_{[M]} \sum_{l=0}^{N}(-1)^l N! \left(\begin{smallmatrix}  d_M\\l \end{smallmatrix}\right)[M] \pmod{(q-1)\sH(A)}\\
               &\equiv 0\pmod{(q-1)\sH(A)}.
              \end{aligned}$$
      
          Now, we will show that $\kappa_1$ is injective. Suppose $\sum_{\lambda\in \aleph}  b_{\lambda}[S_{i_1}]^{l_1}[S_{i_2}]^{l_2}\cdots [S_{i_t}]^{l_t}=0$ (finite sum) in $\overline{C\sH(A)}$, then  $\sum_{\lambda\in \aleph}  b_{\lambda} \prod_{j=1}^t [l_j]_q! F_{\lambda}^M [M]\equiv 0\pmod{(q-1)}$ for finitely many $[M]$. Note that $F_{\lambda}^M\equiv F_{\lambda}^M(1)\pmod{(q-1)}$ and $F_{\lambda}^M(1)=\chi_c(V(\lambda;M))$, it follows that $\sum_{\lambda\in \aleph}  b_{\lambda}\prod_{j=1}^t (l_j!)\chi_c(V(\lambda;M))\equiv 0\pmod{(q-1)}$. Then $\sum_{\lambda\in \aleph}\bar{b_{\lambda}}$ $\bar{f_{i_1}}^{l_1}\bar{f_{i_2}}^{l_2}\cdots \bar{f_{i_t}}^{l_t}=0$ in $\overline{U^+_{\Z}}$, the proof is completed.
        \eproof
        \elemma
      
        According to \cite[Lemma 12.11]{Lusztig2000}, we have the following
        \blemma
          The morphism $\kappa_2: \overline{U^+_{\Z}}\lrw \overline{C\sH(\varLambda_Q)}$ sending $f_i$ to $[S_i]$ is a surjective homomorphism of $\Z/{(q-1)}$-algebras.
        \bproof
          It suffices to show that for $i\neq j\in I$, there are the following relations
                 $$\sum_{l=0}^{N}(-1)^l \left(\begin{smallmatrix}N \\ l\end{smallmatrix}\right) [S_i]^{N-l}[S_j] [S_i]^{l}=0 \pmod{(q-1)\sH(\varLambda_Q)}.$$
          in $\overline{C\sH(\varLambda_Q)}$, where $N=1+a_{ij}+a_{ji}$.
      
          If there are $m:=a_{ij}+a_{ji}$ arrows ${h_r}$ from $i\to j$ and $a_{ij}+a_{ji}$ arrows $\bar{h}_r$ from $j\to i$.  Denote by $a_{M,l}$ the multiple of $[M]$ in $[S_i]^{N-l}[S_j] [S_i]^{l}$ for $M=(V,(x_h)_{h\in \bar{Q}_1})\in \mathrm{rep}^n\varLambda_Q$. Set $V_M:=\sum_{\bar{h}_r} \mathrm{im}\bar{h}_r$ and $W_M:=\bigcap_{h_r}\Ker h_r$, then 
           $a_{M,l}\neq 0$ iff $V_M\subset S_i^{\oplus l} \subset W_M$. If so, then $a_{M,l}=[l]_q! [N-l]_q! [\begin{smallmatrix}  d_M-\bar{d}_M\\l-\bar{d}_M \end{smallmatrix}]_q$, here $d_M=\dim W_M\geq 1$, $\bar{d}_M=\dim V_M$. Note that there is a sequence exact except at the fourth item $V_j^{m}$
              $$0\lrw W_M \hookrightarrow V_i  \stackrel{(x_h)}\lrw V_j^{m}\stackrel{(x_{\bar{h}})}\lrw V_M \lrw 0,$$
          it follows that $\dim W_M-\dim V_M\geq \dim V_i-\dim V_j^{m}=N-m=1$. Then we have that 
              $$\begin{aligned}
              &\sum_{l=0}^{N}(-1)^l \left(\begin{smallmatrix}N \\ l\end{smallmatrix}\right) [S_i]^{N-l}[S_j] [S_i]^{l} \\
              &= \sum_{[M]}\sum_{l=0}^{N}(-1)^l \left(\begin{smallmatrix}N \\ l\end{smallmatrix}\right)[l]_q! [N-l]_q! [\begin{smallmatrix}  d_M-\bar{d}_M\\l-\bar{d}_M \end{smallmatrix}]_q [M]\\
              &\equiv  \sum_{[M]}\sum_{l=0}^{N}(-1)^l N! \left(\begin{smallmatrix}  d_M-\bar{d}_M\\l-\bar{d}_M \end{smallmatrix}\right) [M] \pmod{(q-1)\sH(\varLambda_Q)}\\
              &\equiv 0\pmod{(q-1)\sH(\varLambda_Q)}.
              \end{aligned}$$
        \eproof
        \elemma
      
        Notice that $A$ can be regarded as a quotient algebra of $\varLambda_Q$, then there is an epimorphism $Res: \sH(\varLambda_Q)\to \sH(A)$ according to Proposition \ref{prop1}, which induces an epimorphism $\theta: \sH(\varLambda_Q)/(q-1)\sH(\varLambda_Q)\to \sH(A)/(q-1)\sH(A)$. Furthermore, simple $\varLambda_Q$-modules are precisely simple $A$-modules, then the epimorphism $\theta$ maps $[S_i]$ to $[S_i]$ for $i\in Q_0$. Hence, the following diagram commutes.
           $$\begin{tikzcd}
            \overline{U^+_{\Z}} \arrow[d, "\kappa_1"] \arrow[r, "\kappa_2"]   & \overline{C\sH(\Lambda_Q)}\arrow[dl, "\theta" ] \\
            \overline{C\sH(A)}
          \end{tikzcd}
         $$
      
         \btheorem\label{thm1}
            There is an isomorphism of algebras $\varphi^+ : \ssn^+_0(\varLambda_Q)\lrw \ssn^+_0(A)$.
         \bproof
            By Lemma \ref{lem4}, $\kappa_1=\theta\kappa_2$ is an isomorphism, we obtain that $\kappa_2$ is injective. Since $\kappa_2$ is surjective, $\kappa_2$ is also an isomorphism. Then
            $$\theta=\kappa_1\kappa_2^{-1}: \overline{C\sH(\varLambda_Q)}\to \overline{C\sH(A)}$$
            is an isomorphism. Moreover, $\overline{C\sH(\varLambda_Q)}$ (resp. $\overline{C\sH(A)}$) is the universal enveloping algebra of $\ssn^+_0(\varLambda_Q)$ (resp. $\ssn^+_0(A)$), it follows that the isomorphism $Res: \overline{C\sH(\varLambda_Q)}\to \overline{C\sH(A)}$ gives an isomorphism $\varphi^+ : \ssn^+_0(\varLambda_Q)\lrw \ssn^+_0(A)$.
         \eproof
         \etheorem
         Similarly, we have an isomorphism $\varphi^- : \ssn^-_0(\varLambda_Q)\lrw \ssn^-_0(A)$. Combining Theorem \ref{thm1} with Proposition \ref{prop4.3}, we reach the following 
         \btheorem\label{cor4.9}
          As Lie algebras, $\ssn^+_0(\sR_{\varLambda_Q})$ (resp. $\ssn^-_0(\sR_{\varLambda_Q})$) is isomorphic to $\ssn^+_0(A)$ (resp. $\ssn^-_0(A)$).
         \etheorem
      
         \subsection{The Ringel-Hall Lie algebra of a spherical object}
         Now we consider a certain Lie subalgebra of $\ssg(\sR_{\varLambda_Q})$ supported on one single vertex.  Let $\mathrm{thick}(S_i)$ be the smallest triangulated subcategory of $D^b(\mathrm{rep}^n\varLambda_Q)$ containing the simple module $S_i$ and stable under direct summands. Note $S_i$ is a 2-spherical object and $\mathrm{thick}(S_i)$ is an algebraic triangulated category (for $D^b(\mathrm{rep}^n\varLambda_Q)$ is algebraic triangulated category and its triangulated subcategories are also algebraic from \cite[Lemma 3.5.8]{Yang2019}), hence $\mathrm{thick}(S_i)$ is triangle equivalent to the derived DG category $\sD_{fd}(\Gamma_i)$ of DG $\Gamma_i$-modules generated by the unique simple object $k[t]/(t)$, where $\Gamma_i=k[t]$ with $\mathrm{deg}(t)=-1$, see \cite [Theorem 2.2]{Fu2012}. According to  \cite [Theorem 3.2]{Fu2012}, the orbit category $\sD_{fd}(\Gamma_i)/[2]$ admits a canonical triangle structure,  hence $\sD_{fd}(\Gamma_i)/[2]\simeq \mathrm{thick}(S_i)/[2]$ is a triangulated subcategory of $\sR_{\varLambda_Q}$. Let $\ssL_i$ be the Ringel–Hall Lie algebra of $\mathrm{thick}(S_i)/[2]$, then $\ssL_i$ is a Lie subalgebra of $\ssg(\sR_{\varLambda_Q})$.
         
         Let $\bigtriangleup$ be the cyclic quiver with 2 vertices and $\sT_2$ be the  category of finitely generated nilpotent $k\bigtriangleup$-modules. Then $\sD_{fd}(\Gamma_i)/[2]$ is triangulated equivalent to the cluster tube $D^b(\sT_2)/{\tau\circ [-1]}$ of rank 2 by \cite[Lemma 4.6]{Fu2012}. Hence each indecomposable objects in $\sD_{fd}(\Gamma_i)/[2]$ is of the form $\langle n \rangle$ or $\langle -n \rangle$, where $\langle n \rangle$ (resp. $\langle -n \rangle$) is the unique indecomposable $k\bigtriangleup$-module of length $n$ with socle the simple module corresponding to vertex 1 (resp.2).
         By \cite[Theorem 4.14]{Fu2012}, $\ssL_i$ is the Lie algebra over  $\Z/(q-1)$  with the basis $\{h_{i}, \hat{u}_{i,\langle n\rangle},n\in \Z^{\times }\}$ and structure constants given by (in the following $x, y \in \N, m, n \in \Z^{\times}$ and $h_i:=\frac{h_{S_i}}{d(S_i)}$)
         \begin{itemize}
          \item[(i)] $[\hat{u}_{i,<m>},\hat{u}_{i,\langle n \rangle}]=0$ for $m$ and $n$ even;
          \item[(ii)] $[\hat{u}_{i,<m>},\hat{u}_{i,\langle n \rangle}]=0$ for $m$ and $n$ both odd of the same sign;
          \item[(iii)] $[\hat{u}_{i,<2x>},\hat{u}_{i,\langle 2y-1 \rangle}]$=
           $\begin{cases} \hat{u}_{i,\langle 2(x+y)-1 \rangle}+\hat{u}_{i,\langle 2(y-x)-1 \rangle}, & \mbox{for } x<y, \\  \hat{u}_{i,\langle 2(x+y)-1 \rangle}-\hat{u}_{i,\langle 2(x-y)+1 \rangle}, & \mbox{for } x\geq y; \end{cases}$
           \item[(iv)] $[\hat{u}_{i,<2x>},\hat{u}_{i,\langle -2y+1 \rangle}]$=
           $\begin{cases} -\hat{u}_{i,\langle -2(x+y)+1 \rangle}-\hat{u}_{i,\langle 2(x-y)+1 \rangle}, & \mbox{for } x<y, \\  -\hat{u}_{i,\langle -2(x+y)+1 \rangle}+\hat{u}_{i,\langle 2(y-x)-1 \rangle}, & \mbox{for } x\geq y; \end{cases}$
           \item[(v)] $[\hat{u}_{i,\langle -2x \rangle},\hat{u}_{i,\langle 2y-1 \rangle}]$=
           $\begin{cases} -\hat{u}_{i,\langle 2(x+y)-1 \rangle}-\hat{u}_{i,\langle 2(y-x)-1 \rangle}, & \mbox{for } x<y, \\  -\hat{u}_{i,\langle 2(x+y)-1 \rangle}+\hat{u}_{i,\langle 2(x-y)+1 \rangle}, & \mbox{for } x\geq y; \end{cases}$
           \item[(vi)] $[\hat{u}_{i,\langle -2x \rangle},\hat{u}_{i,\langle -2y+1 \rangle}]$=
           $\begin{cases} \hat{u}_{i,\langle -2(x+y)+1 \rangle}+\hat{u}_{i,\langle 2(x-y)+1 \rangle}, & \mbox{for } x<y ,\\  \hat{u}_{i,\langle -2(x+y)+1 \rangle}-\hat{u}_{i,\langle 2(y-x)-1 \rangle}, & \mbox{for } x\geq y; \end{cases}$
           \item[(vii)] $[\hat{u}_{i,\langle 2x-1 \rangle},\hat{u}_{i,\langle -2y+1 \rangle}]$=
           $\begin{cases} \hat{u}_{i,\langle 2(x+y)-2 \rangle}-\hat{u}_{i,\langle -2(x+y)+2 \rangle}+\hat{u}_{i,\langle 2(x-y) \rangle}-\hat{u}_{i,\langle 2(y-x) \rangle}, & \mbox{for } x<y, \\  -h_i+\hat{u}_{i,<4x-2>}-\hat{u}_{i,<-4x+2>}, & \mbox{for }x=y, \\ \hat{u}_{i,\langle 2(x+y)-2 \rangle}-\hat{u}_{i,\langle -2(x+y)+2 \rangle}+\hat{u}_{i,\langle 2(y-x) \rangle}-\hat{u}_{i,\langle 2(x-y) \rangle}, & \mbox{for } x>y; \end{cases}$
           \item[(viii)] $[h_i,\hat{u}_{i,\langle n\rangle}]$=
           $\begin{cases} 0, & \mbox{for $n$ even}, \\  2\hat{u}_{i,\langle n \rangle}, & \mbox{for $n$ positive odd}, \\ -2\hat{u}_{i,\langle n \rangle}, &\mbox{for $n$ negative odd}.\end{cases}$
         \end{itemize}
          We replace the symmetric Euler form in the last term above by its half, which is precisely the bilinear form $(-|-)_{\sR_{\varLambda_Q}}$ we used.
         \bexample\label{ex4.10}
          Let $\hat{u}_{i,\langle 1 \rangle}$ represent $\hat{u}_{S_i}$ and $\hat{u}_{i,\langle -1 \rangle}$ represent $\hat{u}_{S_i[1]}$. Then
          $$[\hat{u}_{S_i},\hat{u}_{S_i[1]}]=-h_i+\hat{u}_{E_i}-\hat{u}_{E_i[1]}$$
          where $E_i[1]=cone(t_i:S_i\to S_i[2])\in D^b(\mathrm{rep}^n\varLambda_Q)$. Comparing with the identity
              $$[\hat{u}_{i,\langle 1 \rangle},\hat{u}_{i,\langle -1 \rangle}]=-h_i+\hat{u}_{i,\langle 2 \rangle}-\hat{u}_{i,\langle -2 \rangle}.$$
         we know  $\hat{u}_{i,\langle 2 \rangle}=\hat{u}_{E_i}$ and  $\hat{u}_{i,\langle -2 \rangle}=\hat{u}_{E_i[1]}$. Indeed, to compute $[S_i][S_i[1]]$, we consider the triangle in $\sR_{\varLambda_Q}$
                $$S_i[1]\lrw L\lrw S_i\stackrel{f}\lrw S_i[2],$$
        where $f\in \Hom_{D^b{\varLambda_Q}}(S_i,S_i)\oplus \Hom_{D^b{\varLambda_Q}}(S_i,S_i[2])$. Note $\Hom_{\varLambda_Q}(S_i,S_i)\cong k\cong \operatorname{Ext}^2_{\varLambda_Q}(S_i,S_i)$, for any $0\neq f=(f_0,f_2)$, $(f_0,f_2)\sim (1,0)$ under $\mathrm{Aut}_{\sR_{\varLambda_Q}}(S_i)$-action, if $f_0\neq 0$. So $L\cong 0$ and $F_{S_i,S_i[1]}^{0}=1$; if $f_0=0$, then $(0,f_2)\sim (0,1)$ , so $L[1]\cong cone(t_i:S_i\to S_i[2])\in D^b(\mathrm{rep}^n\varLambda_Q)$ and $F_{S_i,S_i[1]}^{L}=1$. Therefore
              $$[S_i][S_i[1]]=1+[E_i].$$
         Similarly, we can show
             $$[S_i[1]][S_i]=1+[E_i[1]].$$
         \eexample
      
          \bremark
          $\ssg_0(\sR_A)$ is not isomorphic to $\ssg_0(\sR_{\varLambda_Q})$ as Lie algebras because the root category $\sR_{\varLambda_Q}$ is not proper, there are more elements belonging to $\ssg_0(\sR_{\varLambda_Q})$. Indeed, the last example \ref{ex4.10} shows that $\hat{u}_{E_i}-\hat{u}_{E_i[1]}=[\hat{u}_{S_i},\hat{u}_{S_i[1]}]+h_i\in \ssg_0(\sR_{\varLambda_Q})$ does not belong to $\ssg_0(\sR_A)$, since $\hat{E_i}$ and $\hat{E_i[1]}$ are zero objects in $K_0(\sR_{\varLambda_Q})$.
         \eremark
         Let $\sI$ be an ideal of $\ssg_0(\sR_{\varLambda_Q})$ generated by $\hat{u}_{E_i}-\hat{u}_{E_i[1]}$, $i\in Q_0$. Let $\ssL_{i,0}$ be the Lie subalgebra of $\ssL_i$ generated by $\hat{u}_{i,\langle \pm 1 \rangle}$ and $h_i$, $\sI_i$ be the ideal of $\ssL_{i,0}$ generated by $\hat{u}_{i,\langle 2 \rangle}-\hat{u}_{i,\langle -2 \rangle}$. Denote by $\bar{u}_{i,\langle \pm 1\rangle}$ the image of $\hat{u}_{i,\langle \pm 1\rangle}$ in $\ssL_{i,0} /\sI_i$. Then we have the following 
      
         \blemma\label{Lem4.13}
          The Lie algebra $\ssL_{i,0}/\sI_i$ is the Lie algebra over $\Z/{(q-1)}$ with a basis
          $$\{\bar{h}_{i}, \bar{u}_{i,\langle 1 \rangle},\bar{u}_{i,\langle -1 \rangle}\}$$
          subject to relations
            $$[\bar{h},\bar{u}_{i,\langle \pm 1 \rangle}]=\pm 2\bar{u}_{i,\langle \pm 1 \rangle}\ \ and\ \ [\bar{u}_{i,\langle 1 \rangle},\bar{u}_{i,\langle -1 \rangle}]=\bar{h}_i$$
          Equivalently, $\ssL_i/\sI_i$ is the Lie algebra $\mathfrak{sl}_2$.
          \bproof
           Compute with given relations of $\ssL_i$, we have
           $$\overline{[\hat{u}_{i,\langle 1 \rangle},\hat{u}_{i,\langle -1 \rangle}]}=\overline{h_i+\hat{u}_{i,\langle 2 \rangle}-\hat{u}_{i,\langle -2 \rangle}}=\bar{h}_i.$$
           and
            $$\overline{[h_i,\hat{u}_{i,\langle \pm 1\rangle}]}=\pm 2 \bar{u}_{i,\langle \pm 1\rangle}.$$
         Therefore there exists an algebra homomorphism $\Psi_i$ from $\mathfrak{sl}_2$ to $\ssL_{i,0}/\sI_i$, which is clearly surjective. On the other hand, using the defining relations (i)-(viii) of $\ssL_i$, we can deduce that $\hat{u}_{i,\langle 1\rangle}\notin \sI_i$. So $\Psi_i$ must be injective.
         \eproof
         \elemma

         It is well known that $\ssg_0(\sR_A)$ is a Kac-Moody Lie algebra over $\Z/(q-1)$. Namely, $\ssg_0(\sR_A)$ is the Lie algebra over $\Z/(q-1)$ with generators $\{ u_{S_i}, u_{S_i[1]},h_i|i\in Q_0\}$ subject to the following relations:
         $$ \begin{aligned}
           &[h_i,h_j]=0,\qquad [h_i, u_{S_j}]= (\hat{S_i},\hat{S_j})_A u_{S_j},\qquad [h_i, u_{S_j[1]}]= -(\hat{S_i},\hat{S_j})_A u_{S_j[1]},\\
           &ad(u_{S_i})^{1+c_{ij}}(u_{S_j})=0,\qquad\qquad  ad(u_{S_i[1]})^{1+c_{ij}}(u_{S_j[1]})=0,\ if\ i\neq j,\\
          \end{aligned}$$
          where $c_{ij}=\dim\operatorname{Ext}^1_A(S_i,S_j)+\dim\operatorname{Ext}^1_A(S_j,S_i)$. Note that we have shown that $\ssn_0^{+}(\sR_{\varLambda_Q})$ is isomorphic to $\ssn^+_0(A)$ as Lie algebras in Theorem \ref{cor4.9}. As a result, we have an algebra homomorphism 
              $$\Psi : \ssg_0(\sR_A)\lrw \ssg_0(\sR_{\varLambda_Q})/\sI.$$
          Since the underlying diagram of the quiver $Q$ is connected, the ideal $\Ker \Psi$ is either $\ssg_0(\sR_A)$ or contained in the center of $\ssg_0(\sR_A)$ by \cite[Proposition 1.7]{Kac1984}.

          Let $\ssh'\subset \ssh$ be the center of $\ssg_0(\sR_A)$, then $\ssh'$ in $\ssg_0(\sR_{\varLambda_Q})$ is also a part of the center of $\ssg_0(\sR_{\varLambda_Q})$. Because $\ssg_0(\sR_A)$ and $\ssg_0(\sR_{\varLambda})$ share same killing form $(-|-)_A=(-|-)_{\sR_{\varLambda_Q}}$.

         \begin{theorem}
          Assume $\ssg_0(\sR_{\varLambda_Q})/\sI \neq 0$, then 
              $$\Psi: \ssg_0(\sR_A)/\ssh' \lrw \ssg_0(\sR_{\varLambda_Q})/(\sI+\ssh'), u_{S_i}\mapsto \hat{u}_{S_i} ,\ u_{S_i[1]}\mapsto \hat{u}_{S_i[1]}\ and\ h_i\mapsto h_i.$$
         \end{theorem}

\section*{Acknowledgement}
   We would like to thank Professor Haicheng Zhang for his valuable suggestions and comments, which greatly makes this paper more readable. The second author thanks Quan Situ for helpful conversations and the idea for proof of Lemma \ref{lem3}.

\bibliographystyle{plain}
\bibliography{myref}
\end{document}